\documentclass{article}

\usepackage{graphicx,amsmath,amssymb,extarrows,algorithm,color,hyperref}
\usepackage[noend]{algpseudocode}

\newtheorem{lemma}{Lemma}
\newtheorem{theorem}{Theorem}
\newtheorem{corollary}{Corollary}
\newcommand{\ep}{\hfill $\square$}

\author{Franti\v sek Kardo\v s \\
LaBRI, University of Bordeaux, France \\
\texttt{frantisek.kardos@labri.fr}}
\title{A computer-assisted proof of Barnette-Goodey conjecture: Not only fullerene graphs are Hamiltonian.}
\begin{document}
\maketitle
\begin{abstract}
Fullerene graphs, i.e., 3-connected planar cubic graphs with pentagonal and hexagonal faces, are conjectured to be Hamiltonian. This is a special case of a conjecture of Barnette and Goodey, stating that 3-connected planar graphs with faces of size at most 6 are Hamiltonian. We prove the conjecture.
\end{abstract}

\section{Introduction}

\begin{figure}[ht]
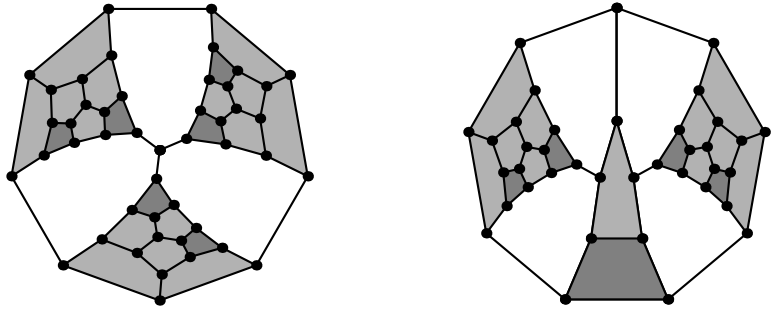

\centerline{
\includegraphics{tutte.1}
\hfil
\includegraphics{tutte.2}
}
\caption{Tutte's first example of a non-Hamiltonian cubic polyhedral graph (left); one of minimal examples on 38 vertices (right).}
\label{fig:tutte}
\end{figure}

Tait conjectured in 1880 that cubic polyhedral graphs (i.e., 3-connected planar cubic graphs) are Hamiltonian. The first counterexample to Tait's conjecture was found by Tutte in 1946; later many others were found, see Figure \ref{fig:tutte}. Had the conjecture been true, it would have implied the Four-Color Theorem. 

However, each known non-Hamiltonian cubic polyhedral graph has at least one face of size 7 or more \cite{ABH,zaks}. It was conjectured that all cubic polyhedral graphs with maximum face size at most 6 are Hamiltonian.  In the literature, the conjecture is usually attributed to Barnette (see, e.g., \cite{mal}), however, Goodey \cite{good} stated it in an informal way as well.

This conjecture covers in particular the class of fullerene graphs, 3-connected cubic planar graphs with pentagonal and hexagonal faces only. Hamiltonicity was verified for all fullerene graphs with up to 176 vertices \cite{ABH}. 
Later on, the conjecture in the general form was verified for all graphs with up to 316 vertices \cite{bgm}.
On the other hand, cubic polyhedral graphs having only faces of sizes 3 and 6 or 4 and 6 are known to be Hamiltonian \cite{good,good2}.

Jendrol$\!$' and Owens proved that the longest cycle of a fullerene graph of order $n$ covers at least $4n/5$ vertices \cite{JO}, the bound was later improved to $5n/6-2/3$ by Kr\' al' et al.~\cite{kral} and to $6n/7+2/7$ by Erman et al.~\cite{erman}.
Maru\v{s}i\v{c} \cite{mar} proved that the fullerene graph obtained from another fullerene graph with an odd number of faces by the so-called leapfrog operation (truncation of the dual; replacing each vertex by a hexagonal face) is Hamiltonian. In fact, Hamiltonian cycle in the derived graph corresponds to a decomposition of the original graph into an induced forest and a stable set. We will use similar technique to prove the conjecture in the general case.

In this paper we prove

\begin{theorem}
Let $G$ be a 3-connected planar cubic graph with faces of size at most 6. Then $G$ is Hamiltonian.
\label{th:main}
\end{theorem}


In the next sections, we reduce the main theorem to Theorem \ref{th:bar} and further to Theorem \ref{th:find} and we introduce terminology and techniques used in the proof of Theorem \ref{th:find}.


\section{Preliminaries}

\subsection{First reduction}
\label{sec:pre}

A \emph{Barnette graph} is a 3-connected planar cubic graph with faces of size at most 6, having no triangles and no two adjacent quadrangles.

We reduce Theorem \ref{th:main} to the case of Barnette graphs:

\begin{theorem}
Let $G$ be a Barnette graph on at least 318 vertices. Then $G$ is Hamiltonian.
\label{th:bar}
\end{theorem}

\begin{lemma}
Theorem \ref{th:bar} implies Theorem \ref{th:main}.
\end{lemma}

Proof. Suppose Theorem \ref{th:bar} true. Let $G$ be a smallest counterexample to Theorem \ref{th:main}. We know that $G$ has at least 318 vertices, since Theorem \ref{th:main} has already been verified for all cubic planar graphs with faces of size at most 6 on at most 316 vertices \cite{bgm}. (The number of vertices of a cubic graph is always even.)

\begin{figure}[ht]
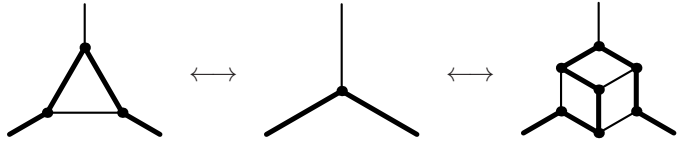

\centerline{
\begin{tabular}{c}
\includegraphics{easy.1}
\end{tabular}
$\longleftrightarrow$
\begin{tabular}{c}
\includegraphics{easy.2}
\end{tabular}
$\longleftrightarrow$
\begin{tabular}{c}
\includegraphics{easy.3}
\end{tabular}
}
\caption{A triangle, as well as three quadrangles sharing a vertex, can be reduced to a single vertex. 
}
\label{fig:easy}
\end{figure}

Assume $f=v_1v_2v_3$ is a triangle in $G$. If one of the faces adjacent to $f$ is a triangle, then, by 3-connectivity, $G$ is (isomorphic to) $K_4$, a Hamiltonian graph.
Therefore, all the three faces adjacent to $f$ are of size at least $4$. Let $G_1$ be a graph obtained from $G$ by replacing $v_1v_2v_3$ by a single vertex $v$. It is easy to see that $G_1$ is a 3-connected cubic planar graph with faces of size at most 6, moreover, every Hamiltonian cycle of $G_1$ can be extended to a Hamiltonian cycle of $G$, see Figure \ref{fig:easy} for illustration.

From this point on we may assume that $G$ contains no triangles.
Let $f_1$ and $f_2$ be two adjacent faces of size 4 in $G$. Let $v_1$ and $v_2$ be the vertices they share; let $f_1=v_1v_2u_3u_4$, let $f_2=v_1v_2w_3w_4$. We denote by $f_3$ (resp. $f_4$) the face incident to $u_3$ and $w_3$ ($u_4$ and $w_4$, respectively). If both $f_3$ and $f_4$ are quadrangles, then, by 3-connectivity, $G$ is the graph of a cube, which is Hamiltonian. Suppose $d(f_4)\ge 5$ and $d(f_3)=4$. Let $G_2$ be a graph obtained from $G$ by collapsing the faces $f_1$, $f_2$, $f_3$ to a single vertex. Again, $G_2$ is a 3-connected cubic planar graph with faces of size at most 6, moreover, every Hamiltonian cycle of $G_2$ can be extended to a Hamiltonian cycle of $G$, see Figure \ref{fig:easy}.

\begin{figure}[ht]
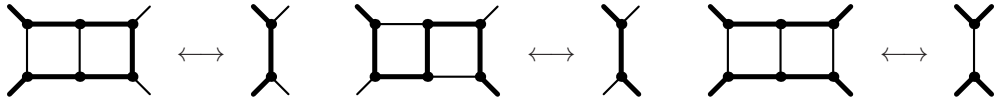

\centerline{
\begin{tabular}{c}
\includegraphics{easy.4}
\end{tabular}
$\longleftrightarrow$
\begin{tabular}{c}
\includegraphics{easy.7}
\end{tabular}
\hfil
\hfil
\begin{tabular}{c}
\includegraphics{easy.5}
\end{tabular}
$\longleftrightarrow$
\begin{tabular}{c}
\includegraphics{easy.8}
\end{tabular}
\hfil
\hfil
\begin{tabular}{c}
\includegraphics{easy.6}
\end{tabular}
$\longleftrightarrow$
\begin{tabular}{c}
\includegraphics{easy.9}
\end{tabular}
}
\caption{A pair of adjacent quadrangles can be reduced to a single edge. 
}
\label{fig:easy2}
\end{figure}

Finally, suppose that both $f_3$ and $f_4$ are of size at most 5. We remove the vertices $v_1$ and $v_2$, identify $u_3$ with $w_3$ and $u_4$ with $w_4$; in this way we obtain a graph $G_3$. It can be verified that $G_3$ is a 3-connected cubic planar graph with all the faces of size at most 6, unless $G$ is the 12-vertex graph obtained from the cube by replacing two adjacent vertices by triangles, which is impossible since $G$ has no triangles. Again, every Hamiltonian cycle of $G_3$ can be extended to a Hamiltonian cycle of $G$, as seen on Figure \ref{fig:easy2}. \ep

\subsection{Cyclic edge-connectivity of Barnette graphs}

Let $G$ be a graph. For a set of vertices $X$, we denote $G[X]$ the subgraph of $G$ induced by $X$.
For a set of vertices $X$, $\emptyset \ne X \ne V(G)$, the set of edges of $G$
having exactly one end-vertex in $X$ form a cut-set of $G$. An edge-cut $(X,Y)$, where $Y=V(G)\setminus X$, is \emph{cyclic} if both $G[X]$ and $G[Y]$ contain a cycle. Finally, a graph is cyclically $k$-edge-connected if it has no cyclic edge-cuts of size smaller than $k$.

\begin{lemma}
Let $G$ be a Barnette graph. Then $G$ is cyclically $4$-edge-connected.
\label{le:c4ec}
\end{lemma}

Proof. Suppose that $G$ contains a cyclic $3$-edge-cut $(X,Y)$. Choose $X$ inclusion-wise minimal. It is easy to see that the cut-edges are pairwise non-adjacent. Let $x_1$, $x_2$, $x_3$ be the vertices of $X$ incident to the cut-edges. We prove that they are pairwise non-adjacent: Suppose that two of them, say $x_1$ and $x_2$, are adjacent. Then, by minimality of $X$, $X^\prime=X\setminus\{x_1,x_2\}$ is acyclic with $(X^\prime,V(G)\setminus X^\prime)$ being a 3-edge-cut, and hence, $|X^\prime|=1$, $X^\prime=\{x_3\}$, so thus $G[X]$ is a triangle, which is impossible in a Barnette graph. 

Let $y_i$ be the other endvertex of the cut-edge incident to $x_i$, $i=1,2,3$. We prove that these three vertices are also pairwise non-adjacent:
Since $G$ has no triangles, $G[\{y_1,y_2,y_3\}]$ has at most two edges. If it had exactly two edges, then $G$ would contain a 2-edge-cut, which is impossible. Suppose now that $y_1$ and $y_2$ are adjacent, but $y_3$ is not adjacent to any of them. Each of the two faces incident to the edge $x_3y_3$ has at least three incident vertices in both $X$ and $Y$, therefore, it is a hexagon, and there are exactly three incident vertices in both $X$ and $Y$. Let $z_i$ be the common neighbor of $y_3$ and $y_i$, $i=1,2$. Then $z_1$ and $z_2$ are adjacent, otherwise there would be a 2-edge-cut in $G$. But then $y_3z_1w_2$ is a triangle in $G$, a contradiction.

As $y_1$, $y_2$, $y_3$ are pairwise non-adjacent, for each face incident to any cut-edge, there are at least three incident vertices in both $X$ and $Y$, therefore, each such face is a hexagon having three incident vertices in both $X$ and $Y$. Let $x_{ij}$ be the common neighbor of $x_i$ and $x_j$, $1\le i < j \le 3$. By minimality of $X$, $X^\prime = X\setminus \{x_1,x_2,x_3,x_{12},x_{13},x_{23}\}$ is a single vertex, and so $G[X]$ is the union of three 4-faces pairwise adjacent to each other, which is impossible in a Barnette graph. \ep

\subsection{Goldberg vectors, Coxeter coordinates, and nanotubes}

Let $f_1$ and $f_2$ be two faces of an infinite hexagonal grid $H$. Then there is a (unique) translation $\phi$ of $H$ that maps $f_1$ to $f_2$. The vector $\vec{u}$ defining $\phi$ can be expressed as an integer combination of two \emph{unit vectors} -- those that define translations mapping a hexagon to an adjacent one. Out of the six possible unit vectors, we choose a pair $\vec{u}_1,\vec{u}_2$ making a $60^\circ$ angle such that $f_2$ is inside this angle starting from $f_1$. Then the coordinates $(c_1,c_2)$ of $\vec{u} = c_1\vec{u}_1 + c_2\vec{u}_2$ are non-negative integers,  called the \emph{Coxeter coordinates} of $\phi$ \cite{cox}. 

We may always assume that $c_1\ge c_2$. The pair $(c_1,c_2)$ determines the mutual position of a pair of hexagons in a hexagonal grid, it is also called a \emph{Goldberg vector}. Observe that, for example, $(1,0)$ corresponds to a pair of adjacent faces, $(1,1)$ corresponds to a pair of non-adjacent faces with an edge connecting them (and thus having two distinct common neighboring faces), whereas $(2,0)$ corresponds to a pair of non-adjacent faces with two paths of length 2 connecting them (and thus sharing a single common neighboring face), etc.

The Coxeter coordinates are used to define nanotubical graphs in the following way:

Let $(c_1,c_2)$ be a pair of integers with $c_1\ge c_2$. Fix a pair of unit vectors $\vec{u}_1$ and $\vec{u}_2$ making a $60^\circ$ angle. A graph obtained from an infinite hexagonal grid by identifying objects (vertices, edges, and faces) whose mutual position is (an integer multiple of) the vector $c_1\vec{u}_1+c_2\vec{u}_2$ is the \emph{infinite nanotube} of \emph{type} $(c_1,c_2)$.

If $c_1+c_2\le 2$ then the infinite nanotube is not 3-connected. Since nanotubes with $c_1+c_2=3$ contain cyclic 3-edge-cuts and Barnette graphs are cyclically 4-edge-connected, we will only be interested in nanotubes with $c_1+c_2\ge 4$.

Let $N$ be an infinite nanotube of type $(c_1,c_2)$. Let $f_1$ and $f_2$ be two hexagons of the hexagonal grid $H$ at mutual position $(c_1,c_2)$ corresponding to the same hexagon $f$ of $N$. Let $P$ be a dual path of length $c_1+c_2$ connecting the vertices $f_1^*$ and $f_2^*$ in $H^*$. Then the edges corresponding to the edges of $P$ form a cyclic edge-cut in $H$ of cardinality $c_1+c_2$. A cyclic sequence of hexagonal faces of $N$ corresponding to the vertices of $P$ is called a \emph{ring} in $N$. 

A finite 2-connected subgraph of an infinite nanotube is an \emph{open-ended nanotube} if it contains at least one ring. A Barnette graph is a \emph{nanotube} if it contains an open-ended nanotube of some type as a subgraph.
Observe that the same graph may be considered as a nanotube of more than one type.

\begin{figure}[ht]
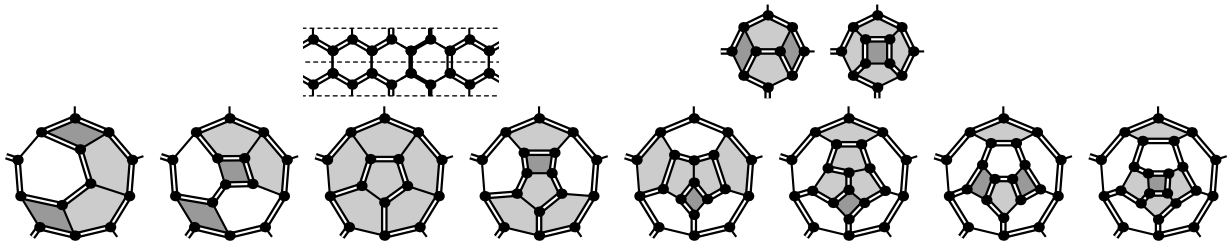

\centerline{
\includegraphics[scale=1.0]{./Figures/caps50.0}
\hfil
\includegraphics[scale=1.0]{./Figures/caps40.1}
\includegraphics[scale=1.0]{./Figures/caps40.2}
}

\centerline{
\includegraphics[scale=1.0]{./Figures/caps50.1}
\includegraphics[scale=1.0]{./Figures/caps50.2}
\includegraphics[scale=1.0]{./Figures/caps50.3}
\includegraphics[scale=1.0]{./Figures/caps50.4}
\includegraphics[scale=1.0]{./Figures/caps50.5}
\includegraphics[scale=1.0]{./Figures/caps50.6}
\includegraphics[scale=1.0]{./Figures/caps50.7}
\includegraphics[scale=1.0]{./Figures/caps50.8}
}
\caption{Two ways to cover the $2p$ vertices separated by two consecutive cyclic $p$-edge-cuts in a $(p,0)$-nanotube by two paths (top left for $p=5$). A path joining two consecutive pending edges covering all the vertices, for every possible cap of $(p,0)$-nanotubes for $p=4$ (top right line) and for $p=5$ (bottom line).}
\label{fig:50tube}
\end{figure}

\begin{figure}[pht]
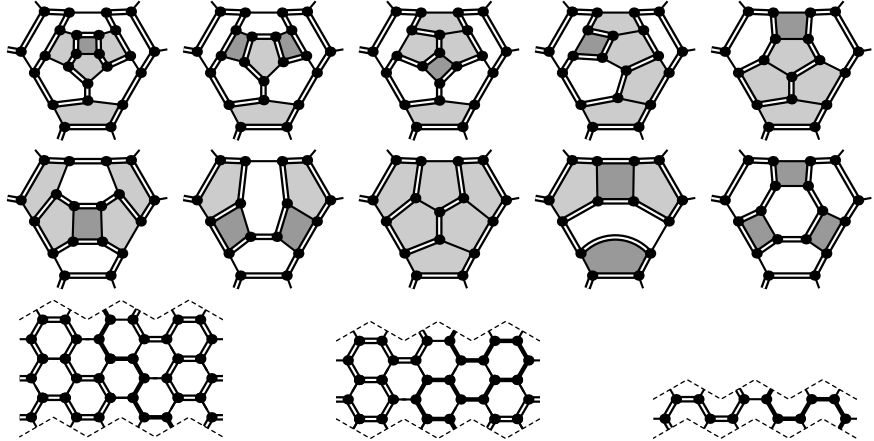

\centerline{
\includegraphics[scale=1.0]{./Figures/caps33.3}
\includegraphics[scale=1.0]{./Figures/caps33.4}
\includegraphics[scale=1.0]{./Figures/caps33.5}
\includegraphics[scale=1.0]{./Figures/caps33.6}
\includegraphics[scale=1.0]{./Figures/caps33.7}
}

\centerline{
\includegraphics[scale=1.0]{./Figures/caps33.8}
\includegraphics[scale=1.0]{./Figures/caps33.9}
\includegraphics[scale=1.0]{./Figures/caps33.10}
\includegraphics[scale=1.0]{./Figures/caps33.11}
\includegraphics[scale=1.0]{./Figures/caps33.12}
}

\centerline{
\includegraphics[scale=1.0]{./Figures/caps33.0}
\hfil
\includegraphics[scale=1.0]{./Figures/caps33.1}
\hfil
\includegraphics[scale=1.0]{./Figures/caps33.2}
}
\caption{For each possible cap of a $(3,3)$-nanotube, a path leaving the cap by a prescribed pair of edges is given (first two rows). For the last cap, we added three hexagons of the tube to make the construction work. To connect the two caps and to cover the tubical part of the graph, it suffices to combine an appropriate number of the first two patterns of the last row (and/or their mirror images) and finish by the third one.}
\label{fig:33tube}
\end{figure}

\begin{figure}[pht]
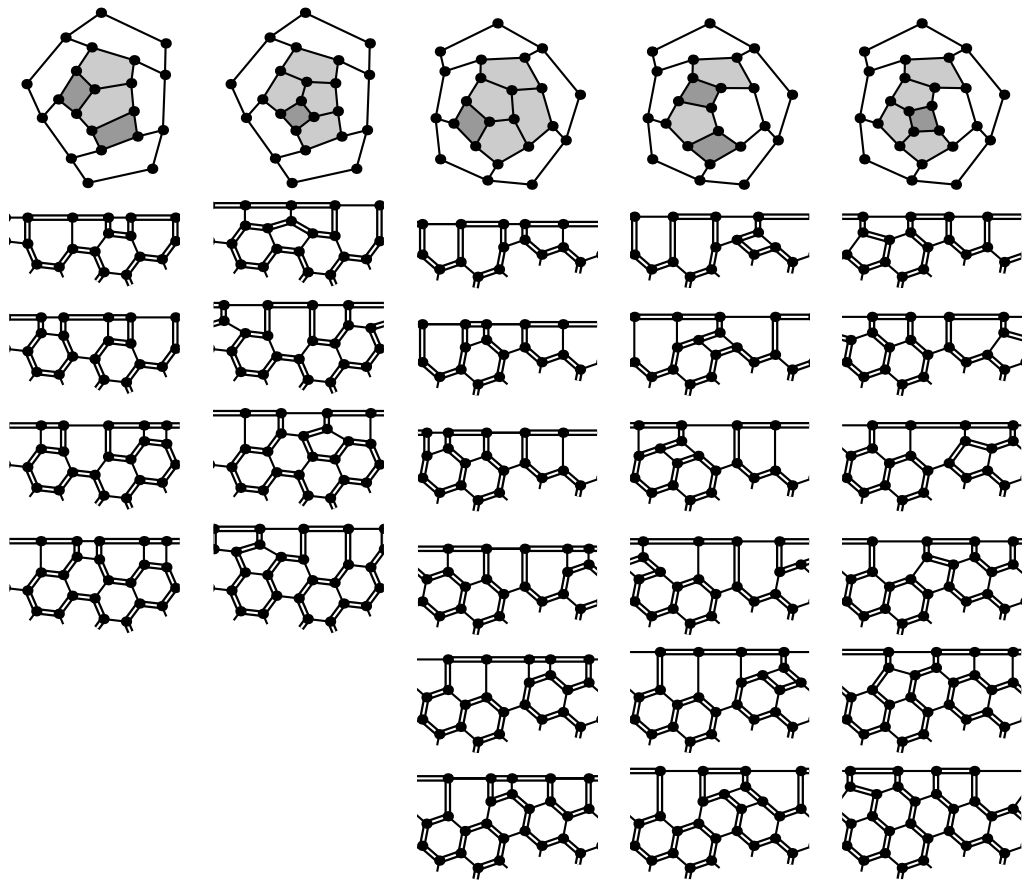

\centerline{
\begin{tabular}{ccccc}
{\includegraphics[scale=1.0]{./Figures/caps32.1}} &{\includegraphics[scale=1.0]{./Figures/caps32.2}} &{\includegraphics[scale=1.0]{./Figures/caps41.1}} &{\includegraphics[scale=1.0]{./Figures/caps41.2}} &{\includegraphics[scale=1.0]{./Figures/caps41.3}}\\
{\includegraphics[scale=1.0]{./Figures/caps32.11}} &{\includegraphics[scale=1.0]{./Figures/caps32.21}} &{\includegraphics[scale=1.0]{./Figures/caps41A.1}} &{\includegraphics[scale=1.0]{./Figures/caps41B.1}} &{\includegraphics[scale=1.0]{./Figures/caps41C.1}}\\
{\includegraphics[scale=1.0]{./Figures/caps32.12}} &{\includegraphics[scale=1.0]{./Figures/caps32.22}} &{\includegraphics[scale=1.0]{./Figures/caps41A.2}} &{\includegraphics[scale=1.0]{./Figures/caps41B.2}} &{\includegraphics[scale=1.0]{./Figures/caps41C.2}}\\
{\includegraphics[scale=1.0]{./Figures/caps32.13}} &{\includegraphics[scale=1.0]{./Figures/caps32.23}} &{\includegraphics[scale=1.0]{./Figures/caps41A.3}} &{\includegraphics[scale=1.0]{./Figures/caps41B.3}} &{\includegraphics[scale=1.0]{./Figures/caps41C.3}}\\
{\includegraphics[scale=1.0]{./Figures/caps32.14}} &{\includegraphics[scale=1.0]{./Figures/caps32.24}} 
&{\includegraphics[scale=1.0]{./Figures/caps41A.4}} &{\includegraphics[scale=1.0]{./Figures/caps41B.4}} &{\includegraphics[scale=1.0]{./Figures/caps41C.4}}\\
&&{\includegraphics[scale=1.0]{./Figures/caps41A.5}} &{\includegraphics[scale=1.0]{./Figures/caps41B.5}} &{\includegraphics[scale=1.0]{./Figures/caps41C.5}}\\
&&{\includegraphics[scale=1.0]{./Figures/caps41A.6}} &{\includegraphics[scale=1.0]{./Figures/caps41B.6}} &{\includegraphics[scale=1.0]{./Figures/caps41C.6}}\\
\end{tabular}
}
\caption{For every cap of a nanotube of types $(3,2)$ (first two columns) and $(4,1)$ (the rest), and for every position of the cap relative to the two paths covering the tubical part of the graph, a completion of the Hamilton cycle in the cap is given. In the first row, the caps are drawn together with the first ring of the tube.}
\label{fig:caps}
\end{figure}

Let $G$ be a nanotube. We call a \emph{cap} any of the two inclusion-wise minimal 2-connected subgraphs of $G$ that can be obtained as a component of a cyclic edge-cut defined by a set of edges intersecting a line perpendicular to the vector defining the corresponding open-ended nanotube. See Figures \ref{fig:50tube}, \ref{fig:33tube}, and \ref{fig:caps} for illustration.

\begin{lemma}
Let $G$ be a Barnette graph which is a nanotube of type $(p_1,p_2)$ with $p_1+p_2=4$. Then $(p_1,p_2)=(4,0)$.
\end{lemma}

We omit the details of the proof, as it is similar to the proof of Lemma \ref{le:c4ec}: It suffices to prove that every (potential) cap of a nanotube of type $(3,1)$ or $(2,2)$ contains a triangle or a pair of adjacent quadrangles.

\begin{lemma}
Let $G$ be a Barnette graph which is a nanotube of type $(p_1,p_2)$ with $(p_1,p_2)\in \{(4,0), (5,0),(4,1),(5,1),(3,2),(4,2),(3,3),(4,3)\}$. Then $G$ is Hamiltonian.
\label{le:smalltubes}
\end{lemma}

Proof. We may suppose that $G$ has at least $318$ vertices (at least 161 faces). Since the caps of the tube are of bounded size (at most 5, 10, 6, 11, 5, 10, 10, 14 faces, respectively, each), the tubical part of $G$ contains a large number of disjoint rings.

We provide a construction of a Hamilton cycle in such graphs: First, we find a pair of paths covering the vertices of the tubical part of $G$; then, we verify that for each possible cap it is always possible to connect the two paths in a way that all the vertices of the cap are covered as well.

In a nanotube of type $(p,0)$, $p\ge 4$, for each $p$-edge-cut corresponding to a ring, we construct the two paths tranversing the tube in a way that each path contains one cut-edge incident to the same hexagonal face. Let us call this hexagon a \emph{transition face}. 
For two consecutive rings, the transition faces are adjacent and once the transition face is fixed for one ring, we are free to choose any of the two adjacent hexagons in the next one to be the transition face, see Figure \ref{fig:50tube} for illustration.

To complete the proof for $(4,0)$- and for $(5,0)$-nanotubes, it suffices to verify that for every possible cap, there exists a path covering all the vertices of the cap leaving the cap by two edges adjacent to the same hexagonal face of the first ring of the tube. Since the tubical part of $G$ is sufficiently long, we can choose a transition face in the first and the last ring of hexagons regardless of the relative position of the two caps.

For nanotubes of type $(3,3)$, the construction is described in Figure \ref{fig:33tube}.

\begin{figure}[pht]
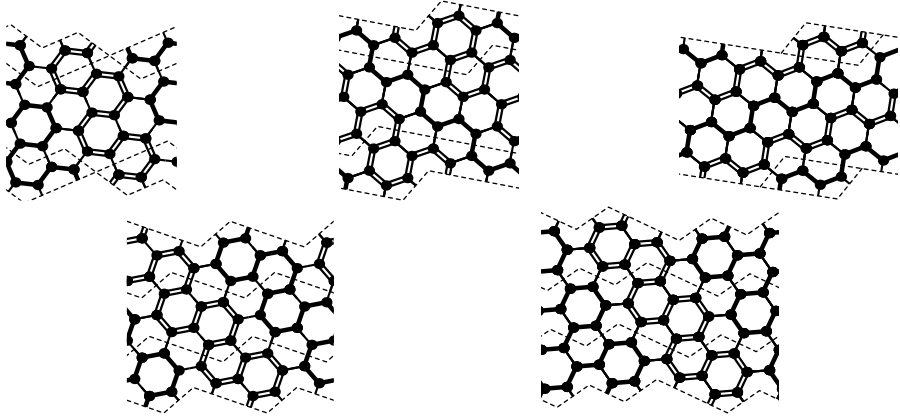

\centerline{
\includegraphics[scale=1.0]{./Figures/caps32.0}
\hfil
\hfil
\includegraphics[scale=1.0]{./Figures/caps41.0}
\hfil
\hfil
\includegraphics[scale=1.0]{./Figures/caps51.0}
}

\centerline{
\includegraphics[scale=1.0]{./Figures/caps42.0}
\hfil
\includegraphics[scale=1.0]{./Figures/caps43.0}
}
\caption{Two paths covering all the vertices of a (potentially infinite) open-ended nanotube of type $(3,2)$, $(4,1)$, $(5,1)$, $(4,2)$, and $(4,3)$, respectively. For each end of the tube, the two dashed lines separate the smallest period of the covering.}
\label{fig:tubes}
\end{figure}

For nanotubes of type $(p_1,p_2)$ with $p_1>p_2>0$, we provide a repetitive pattern to cover the tubical part (see Figure \ref{fig:tubes}) and, for every cap and for every position of the cap with respect to the pattern, a path covering the vertices of the cap (see Figure \ref{fig:caps} for the first two types of nanotubes; we omit the details for the remaining three types). \hfill $\square$

\subsection{Second reduction}

Let $H$ be a plane cubic graph. We denote $H^\parallel $ the $6$-regular multigraph obtained from $H$ by replacing each edge by a pair of parallel edges, equipped with the following black-and-white face-coloring: We color the 2-gons between pairs of parallel edges white and we color the faces of $H^\parallel$ corresponding to the faces of $H$ black. It is easy to see that this is a proper face-coloring of $H^\parallel$.

Let $G$ be a Barnette graph and let $M$ be a perfect matching of $G$. Then $F=E(G)\setminus M$ is a $2$-factor of $G$. A hexagonal face of $G$ incident to three edges of $M$ is called \emph{resonant}.

There is a canonical face-coloring of $G$ with two colors, say black and white, such that each edge of $F$ is incident to one black and one white face. Let $h$ be a white resonant hexagon. Since it is incident to three edges from $M$, the colors of its neighboring faces are alternating black and white.

We transform $F$ into a $6$-regular plane pseudograph in the following way:
First, inside each white resonant hexagon $h$ we introduce a new vertex $v_h$. We remove the three edges incident to $h$ from $F$ and we replace them by six new edges, joining $v_h$ to all the six vertices incident to $h$. Each of the newly created triangles receives the color of the corresponding face adjacent to $h$.
This way we obtain a black-and-white face-colored plane graph with two types of vertices: vertices of degree 2 are the vertices of the underlying Barnette graph, vertices of degree 6 correspond to white resonant hexagons.

\begin{figure}[ht]
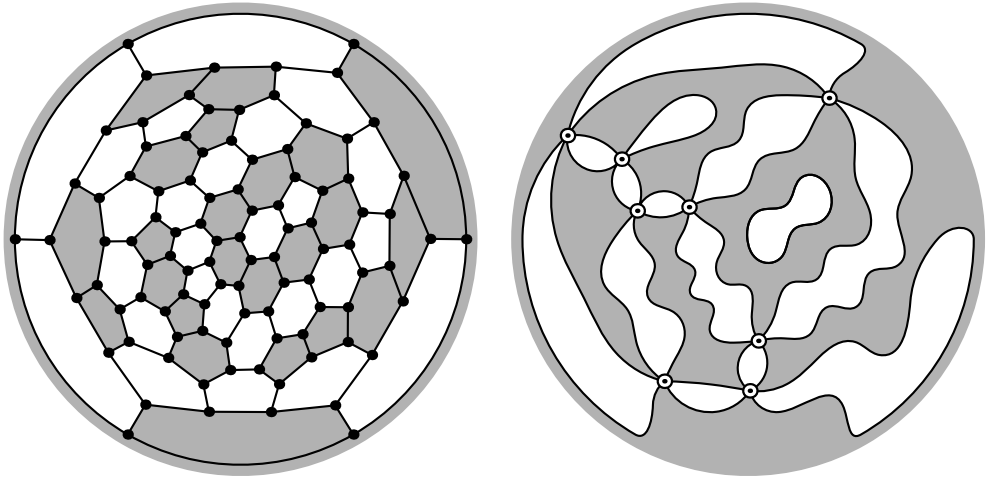

\centerline{
\includegraphics[scale=1]{example.100}\hfil
\includegraphics[scale=1]{example.101}
}
\caption{An example of a black-and-white face-colored 6-regular pseudograph (right) corresponding to a $2$-factor of a Barnette graph (left).}
\label{fig:h6}
\end{figure}

Finally, we suppress all vertices of degree 2. This operation may create loops, parallel edges, and even circular edges incident to no vertex, see Figure \ref{fig:h6} for illustration.
Let $G^M$ be the resulting black-and-white face-colored plane 6-regular pseudograph. 

A $2$-factor $F$ is called $\emph{odd}$ if it consists of an even number of (disjoint) cycles; otherwise it is \emph{even}. The same applies to the corresponding perfect matching.

A $2$-factor $F$ (as well as the corresponding perfect matching $M = E(G)\setminus F$) is called \emph{simple} if $G^M$ has no circular edges and $G^M\cong H^\parallel$ for some cubic planar graph $H$. If this is the case, $H$ is called the \emph{residual graph}.

\begin{lemma}
Let $F$ be a simple $2$-factor of a Barnette graph $G$. Let $n$ be the number of vertices of the corresponding residual graph $H$. If $F$ is odd, then $n=4k+2$ for some $k\ge 1$; otherwise $n=4k$ for some $k\ge 1$.
\end{lemma}

Proof. The number of vertices of a residual graph is always even, since it is a cubic graph. Moreover, the number of cycles in $F$, say $c$, is equal to the number of faces of the residual graph. By Euler's formula,
$$c = 2+|E(H)| - |V(H)| = 2+\frac{3n}2-n = \frac{n+4}{2},$$
so the claim follows immediately.

We will make use of the following classical result:
\begin{theorem}[Payan and Sakarovitch \cite{PS}]
Let $H$ be a cubic graph on $n=4k+2$ vertices ($k\ge 1$). If $H$ is cyclically $4$-edge-connected, then $V(H)$ admits a partition into two sets, say $B$ and $W$, such that $H[B]$ is a stable set and $H[W]$ is a tree.
\label{th:pyber}
\end{theorem}

Observe (by double-counting white-white and black-white edges) that the divisibility condition is a necessary condition for such a partition to exist.
That's why we will only be interested in odd $2$-factors.

\begin{lemma}
Let $G$ be a Barnette graph and let $M$ be an odd simple perfect matching of $G$. If the residual graph is a cyclically $4$-edge-connected, then $G$ is Hamiltonian.
\end{lemma}

Proof. Let $H$ be a cyclically $4$-edge-connected cubic planar graph on $4k+2$ vertices ($k\ge 1$) such that $G^M=H^\parallel$. Let $F=E(G)\setminus M$.  
Recall that vertices of $H$ correspond to white resonant hexagons in $G$ with respect to a fixed cannonical face-coloring of $F$. 
Let $(B,W)$ be a partition of $V(H)$ into an induced (black) stable set $B$ and  an induced (white) tree $W$ given by Theorem \ref{th:pyber}. 

\begin{figure}[pht]
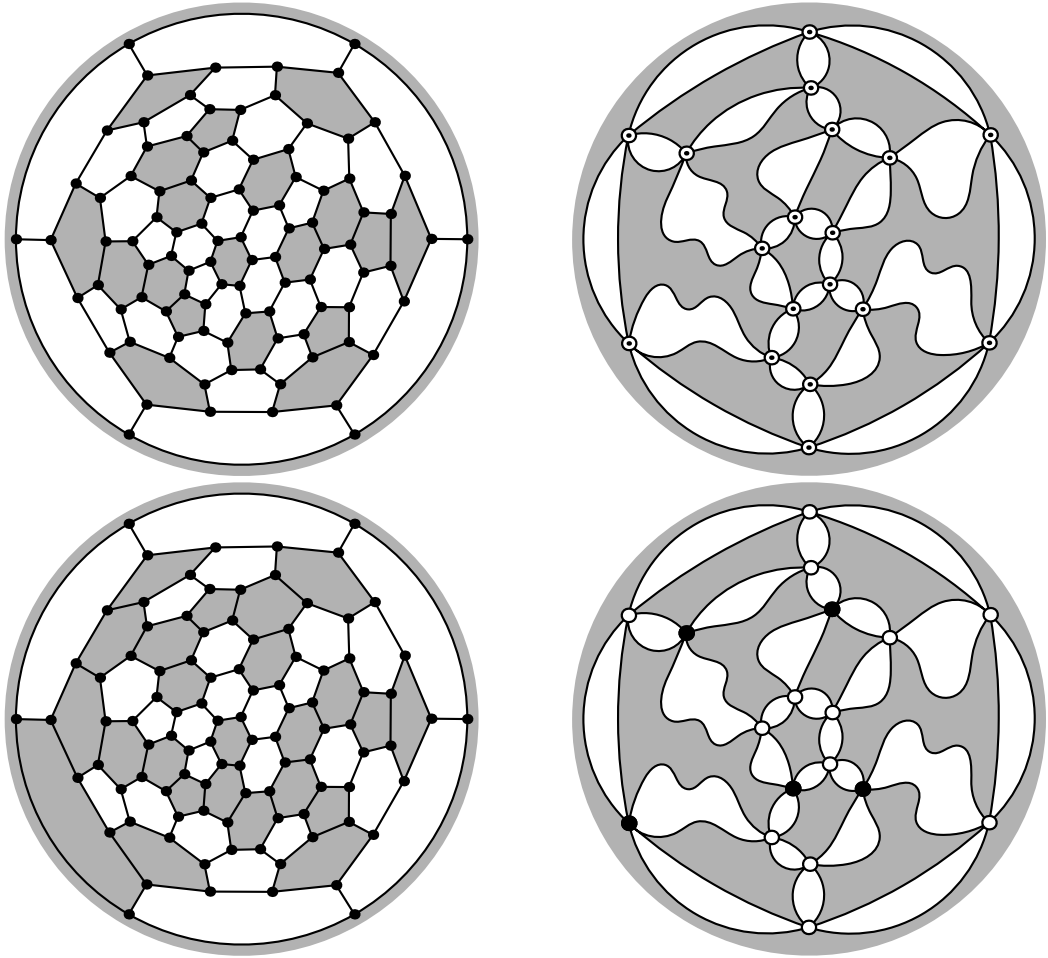

\centerline{
\includegraphics[scale=1]{example.102}\hfill
\includegraphics[scale=1]{example.103}
}
\centerline{
\includegraphics[scale=1]{example.104}\hfill
\includegraphics[scale=1]{example.105}
}
\caption{Clockwise, starting from upper left: An example of a simple $2$-factor $F$ of a Barnette graph $G$; the corresponding planar 6-regular pseudograph $G^M$ which is a double of a cyclically 4-edge connected cubic graph $H$; A decoposition of $H$ into a (black) stable set and a (white) induced tree; the corresponding Hamilton cycle in $G$.}
\label{fig:howItWorks}
\end{figure}

We transform the $2$-factor $F$ and the black-and-white face-coloring of $G$ in the following way:
For each resonant hexagon $h$ corresponding to a black vertex $b$ of $H$, replace the three edges from $F$ incident to $h$ in $G$ by the other three edges; recolor the hexagon $h$ black. Since $B$ induces a stable set in $H$, this operation can be carried out independently for all black vertices of $H$ at once. For each such vertex, the number of edges from $F$ incident to any vertex of $G$ remains unchanged, therefore, $F$ becomes a $2$-factor of $G$, say $F^\prime$.

We claim that it consists of a single cycle. To prove that, it suffices to observe that the graph $(V(G),F^\prime)$ has a single white face (as $H[W]$ is connected) and a single black face (as $H[W]$ is acyclic). See Figure \ref{fig:howItWorks} for illustration. \ep

It remains to prove that such a situation occurs for at least one perfect matching for any Barnette graph not known to be Hamiltonian yet.

\begin{theorem}
Let $G$ be a Barnette graph on at least $318$ vertices. Then there exists an odd simple perfect matching $M$ of $G$ such that the residual graph 
$H$ is cyclically $4$-edge-connected, unless $G$ is a nanotube of type $(4,0)$, $(5,0)$, $(4,1)$, $(5,1)$, $(3,2)$, $(4,2)$, $(3,3)$, or $(4,3)$.
\label{th:find}
\end{theorem}

In the rest of the paper, we prove Theorem \ref{th:find}. 
We describe the general approach in Section \ref{sec:pro}, and we specify the computer-assisted part in Section \ref{sec:com}.

We claim (without proof) that in order to prove Theorem \ref{th:find} it suffices to consider a simple odd $2$-factor maximizing the number of white resonant hexagons.

\subsection{Generalized 2-factors}

We will call a \emph{$2^*$-factor} of a Barnette graph $G$ any spanning subgraph $F$ of $G$ such that 
each component of $F$ is a connected regular graph of degree 1 or 2 -- an isolated edge or a cycle. For a 2$^*$-factor $F$ of a Barnette graph $G$, let $F^{(0)}$ be the set of isolated edges of $F$; let $G^{(2)}$ be a plane graph obtained from $G$ by replacing each edge of $F^{(0)}$ by a 2-gon; let $F^{(2)}$ be the set of edges of $G^{(2)}$ corresponding to those from $F$. Then $F^{(2)}$ is a $2$-factor of $G^{(2)}$ in the common (strict) sense.

Given a 2$^*$-factor $F$ of a Barnette graph $G$, there are two cannonical black-and-white face-colorings of $G^{(2)}$ (complementary to each other) with the following property: an edge $e$ of $G^{(2)}$ is incident to a white and a black face if and only if $e$ belongs to $F^{(2)}$ (otherwise $e$ is incident to two faces of the same color).

A $2^*$-factor $F$ of a Barnette graph $G$ is called \emph{quite good} if for each of the two canonical black-and-white face-colorings of $G^{(2)}$ induced by $F^{(2)}$ the 2-gons corresponding to the edges of $F^{(0)}$ have all the same color. Given a quite good $2^*$-factor of a Barnette graph $G$, we will always assume that a canonical coloring of $G^{(2)}$ such that all the 2-gons of $G^{(2)}$ are black is given along. 

A quite good $2^*$-factor $F$ of a Barnette graph $G$ is called \emph{good} if, after having fixed a planar embedding of $G$ such that the outer face is a white one, no cycle of $F$ is inside another.

Observe that given a good $2^*$-factor $F$ of $G$, for any planar embedding of $G$ with a white outer face, the set of faces inside a fixed cycle $C$ of $F$ is always the same and these faces correspond to a sub-tree of the dual graph $G^*$ (empty if $C$ is a 2-cycle).


\begin{lemma}
Let $F$ be a good $2^*$-factor of a Barnette graph $G$. 
Let $f$ be the number of all the faces of $G$, let $q_k$ be the number of non-resonant white faces of size $k$ in $G$ ($k=4,5,6$); let $c$ be the number of components of $F$. Then $q_5$ is even, moreover, $f+q_4+q_5/2+c \equiv 0 \pmod{2}$. 
\label{l:paths}
\end{lemma}

Proof. Let $n$ be the number of vertices of $G$, let $f_k$ be the number of all faces of size $k$ in $G$, let $x_k$ be the number of black faces of size $k$ in $G$. Euler's formula yields  $n=8+f_5+2f_6$. 
If a cycle covers $c_4\ge 0$ quadrangles, $c_5\ge 0$ pentagons, and $c_6\ge 0$ hexagons, its length is $2+2c_4+3c_5+4c_6$. 

Clearly, each vertex is covered by exactly one cycle, thus we have
$$
8+f_5+2f_6=n=2c+2x_4+3x_5+4x_6 = 2c+2(f_4-q_4)+3(f_5-q_5)+4x_6,
$$ 
since only hexagons can be resonant, and thus $f_k=x_k+q_k$ for $k=4,5$. Therefore, 
$$
8+q_5+2f_6=2c+2(f_4-q_4)+2(f_5-q_5)+4x_6,
$$ so $q_5$ is even. By dividing by two and rearranging the terms we obtain 
$$
4+f_4+f_5+f_6+q_4+q_5/2 +c= 2c+2f_4+2f_5-q_5+2x_6,
$$
the claim immediately follows. \hfill $\square$

\bigskip

Let $F$ be a good $2^*$-factor in a Barnette graph $G$.
Let us consider the structure of the graph $G^{(2)}$. We introduce an auxiliary graph $\Gamma = \Gamma_G(F)$, defined in the following way: $V(\Gamma)$ is the set of the white non-resonant faces of $G$ (as of $G^{(2)}$). The edges of $\Gamma$ are defined in the next two paragraphs.

Let $C$ be the facial cycle of a (black) 2-gon $f_0$ in $G^{(2)}$. Let $f_0$ be incident to vertices $u$ and $v$ and adjacent to two (white) faces $f$ and $f^\prime$. Then each of $u$ and $v$ is incident to one more face (which has to be white), say $f_u$ and $f_v$, respectively. Since $f_0$ only shares a vertex with $f_u$ and with $f_v$, the faces $f$ and $f^\prime$ are two consecutive white neighbors of $f_u$ ($f_v$). Therefore, the faces $f_u$ and $f_v$ cannot be resonant. We add the edge $f_uf_v$ to $E(\Gamma)$; we call this type of edge of $\Gamma$ \emph{white}.

Let $C$ be a cycle of $F$ (and of $F^{(2)}$) which is not a facial cycle of a face of $G^{(2)}$. It means that $C$ is a boundary of a union of at least two faces of $G$. We consider every pair of adjacent faces inside $C$. Let $f$ and $f^\prime$ be such a pair of faces. Let $u$ and $v$ be the endvertices of the edge incident to both $f$ and $f^\prime$. Then each of $u$ and $v$ is incident to a third face (which has to be white), say $f_u$ and $f_v$, respectively. The faces $f$ and $f^\prime$ are two consecutive black neighbors of $f_u$ ($f_v$). Therefore, the faces $f_u$ and $f_v$ cannot be resonant. We add the edge $f_uf_v$ to $E(\Gamma)$; we call this type of edge of $\Gamma$ \emph{black}.

Observe that for each edge of $\Gamma$, its endvertices are two faces of $G$ at mutual position $(1,1)$. Each edge of $\Gamma$ covers two vertices of $G$ and these pairs of vertices are pairwise disjoint. Therefore, $\Gamma$ is a planar graph. 

Let $f$ be a white pentagon of $G^{(2)}$. It cannot be resonant, so $f$ is a vertex of $\Gamma$. Let $f_1,\dots,f_5$ be the faces adjacent to $f$ (sharing an edge with $f$) in $G^{(2)}$. (Observe that some $f_i$ can be a 2-face: if it is the case, then there is another face $f_i^\prime$ adjacent to $f$ in $G$, and adjacent to $f_i$ in $G^{(2)}$.) Since the size of $f$ is odd, the number of pairs $(f_i,f_{i+1})$ (with $f_6=f_1$) of the same color (both black or both white) has to be odd. If both $f_i$ and $f_{i+1}$ are black, then none of them can be a 2-face, and thus there is a black edge incident to $f$ in $\Gamma$. If both $f_i$ and $f_{i+1}$ are white, then again none of them can be a 2-face, and the vertex incident to $f$, $f_i$, and $f_{i+1}$ is (in $G^{(2)}$) covered by a 2-cycle adjacent both to $f_i$ and $f_{i+1}$, and thus there is a white edge incident to $f$ in $\Gamma$. Altogehter, $f$ is a vertex of odd degree in $\Gamma$.

Similarly, for each non-resonant white hexagon $f$, there is an even number of pairs of consecutive adjacent faces of the same color, hence $f$ is a vertex of non-zero even degree in $\Gamma$. 

A white quadrangle $f$ is always considered non-resonant. Its degree in $\Gamma$ is also always even, however, it can be equal to 0 if the neighboring faces are colored alternatively black and white.

As a result of these local observations, the graph $\Gamma$ can always be edge-decomposed into a set of paths with endvertices at the white pentagons of $G$, a set of cycles, and, eventually, a set of isolated vertices (corresponding to white quadrangles). The number of paths in the decomposition is equal to $q_5/2$, where $q_5$ is the number of white pentagons.

\subsection{Structure of Barnette graphs}

\begin{figure}[ht]
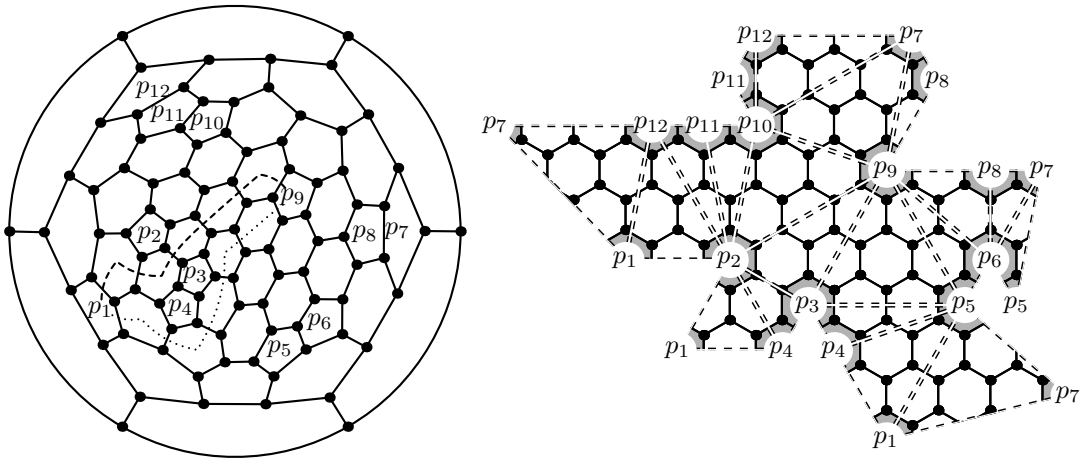

\centerline{
\includegraphics{example.0}
\hfill
\includegraphics{embgrid.1}
}
\caption{An example of a Barnette graph (left). Pentagonal faces are denoted $p_1,\dots,p_{12}$. The mutual position of $p_1$ and $p_9$ is characterized by vectors $(3,3)$ (dotted line) or $(4,2)$ (dashed line). The same graph embedded into a hexagonal grid after being cut along a spanning tree of a triangulation capturing the mutual position of all the small faces (right).}
\label{fig:ex}
\end{figure}

Let $G$ be a Barnette graph and let $p_1$ and $p_2$ be two small faces of $G$. Suppose that there exists an induced dual path $P^*$ connecting $p_1$ and $p_2$ passing only through hexagons. Then if we consider only faces of $G$ corresponding to $P^*$, and if we replace the two small faces by hexagons, we obtain a graph with a cannonical embedding into an infinite hexagonal grid. The Goldberg vector $(c_1,c_2)$ joining the first and the last hexagon is uniquely determined. We will use this vector to characterize the mutual position of $p_1$ and $p_2$ in $G$. 
Observe that the vector of two small faces may depend on the choice of the path joining them, see Figure \ref{fig:ex} for illustration.

\bigskip

Graver \cite{Graver} used the Coxeter coordinates to describe the structure of fullerene graphs. His technique may be extended to a full description of Barnette graphs as well in the following way: A given Barnette/fullerene graph $G$ is represented by a planar triangulation $T$, whose vertices represent the small faces of $G$, and each edge $uv$ is labelled with a Goldberg vector representing the mutual position of the faces represented by $u$ and $v$. The angle between face-adjacent edges (incident to the same triangle of $T$) is well defined and is determined by the labels of the three edges forming the triangle. For a vertex of $T$ representing a pentagon (a quadrangle) the angles around it sum up to $5/3\pi = 300^\circ$ ($4/3\pi = 240^\circ$, respectively).

The existence of a triangulation $T$ is guaranteed by a structural theorem of Alexandrov (see e.g.~\cite{DR}, Theorem 23.3.1, or \cite{pak}, Theorem 37.1), which states (in a more general setting) that any Barnette graph can be embedded onto the surface of a convex (possibly degenerate) polyhedron
so that every face is isometric to a regular polygon with unit edge length;  it suffice then to triangulate the faces of this polyhedron.
Any spanning tree of $T$ may be used to cut the graph $G$ in order to obtain a graph embeddable into the infinite hexagonal grid, see Figure \ref{fig:ex} for illustration. 

We say that a Goldberg vector $\vec{u}=(c_1,c_2)$ is \emph{shorter} than $\vec{u}^\prime=(c_1^\prime, c_2^\prime)$ if and only if the Euclidean length of a segment determined by $\vec{u}$ is shorter than the Euclidean length of a segment determined by $\vec{u}^\prime$ when both embedded into the same hexagonal grid.

Observe that the triangulation representing a Barnette graph is not unique: wherever two adjacent triangles form a convex quadrilateral (once embedded into the hexagonal grid), we may choose the other diagonal of the quadrilateral instead of the existing one as an edge of the triangulation. For example, in the graph depicted in Figure \ref{fig:ex} we could have chosen the edge $p_3p_{10}$ instead of the edge $p_2p_9$, etc.

However, for a triangulation $T$ representing a Barnette graph $G$, the operation switching the diagonals of a convex quadrilateral eventually leads to a triangulation minimal with respect to the sum of lengths of its edges. For example, the triangulation depicted in Figure \ref{fig:ex} is already minimal.

\begin{lemma}
Let $G$ be a Barnette graph, let $T$ be a minimal triangulation representing $G$. Then $T$ has a Hamiltonian path. 
\end{lemma}

Proof. 
Suppose that $T$ has no Hamiltonian path. Then there exists a set $X$ of vertices such that $T\setminus X$ has at least $|X|+2$ connected components. Since $T$ has at most 12 vertices, $|X|\le 5$.

For each component $C$, the set of vertices in $G\setminus C$ having a neighbor in $C$ contains a cycle in $T$ (as $T$ is a triangulation). Therefore, $G\setminus X$ is a plane graph with $|X|\le 5$ vertices and $|V(T)\setminus X|\ge 7$ faces. 




However, a planar graph on at most 5 vertices can have at most 6 faces. (Adding edges increases the number of faces, and (the) planar triangulation on 5 vertices (the triangular bipyramid) only has 6 faces.) 
\ep

Note that the smallest planar graph with desired properties is a bipyramid over a square (which has 6 vertices and 8 faces).

\begin{lemma}
Let $G$ be a Barnette graph, let $T$ be a minimal triangulation representing $G$. Then either $T$ is Hamiltonian, or $T$ can be transformed to a Hamiltonian trangulation by a single diagonal switch.
\label{le:triang}
\end{lemma}

Proof. 
Suppose that $T$ has no Hamiltonian cycle. Then there exists a set $X$ of vertices such that $T\setminus X$ has at least $|X|+1$ connected components. Since $T$ has at most 12 vertices, $|X|\le 5$.

For each component $C$, the set of vertices in $G\setminus C$ having a neighbor in $C$ form a cycle in $T$ (as $T$ is a triangulation). Therefore, $G\setminus X$ is a plane graph with $|X|\le 5$ vertices and $|V(T)\setminus X|\ge 6$ faces. 

There is only one such graph: the triangular bipyramid $B$, which has $5$ vertices and $6$ triangular faces. Out of the six components of $T\setminus B$, at least five are singletons, the sixth may eventually be an isolated edge. It means $T$ has five vertices of degree at least 6, six vertices of degree 3, and eventually a vertex of degree 4.

Let $e=uv$ be an edge of $B$. It is incident to two triangles, each incident to a different component of $T\setminus B$. Let $x$ and $y$ be the vertices of $T\setminus B$ such that $uvx$ and $uvy$ are triangles of $T$. If the quadrilateral $uxvy$ is convex, then the triangulation $T^\prime$ obtained from $T$ by switching $uv$ to $xy$ has at most five vertices of degree 3, so $T^\prime$ has to be Hamiltonian.

It remains to consider the case when for each edge $e$ of $B$, the union of the two incident triangles is a non-convex quadrilateral, meaning that at one of its endvertices, the sum of the angles in the incident triangles is greater than $180^\circ$. Since $B$ has five vertices and nine edges, there is at least one vertex of $B$ with two (disjoint) pairs of incident triangles whose union gives a non-convex angle. But then the sum of the angles around this vertex is greater than $360^\circ$, a contradiction. 
\ep

\begin{figure}[pht]
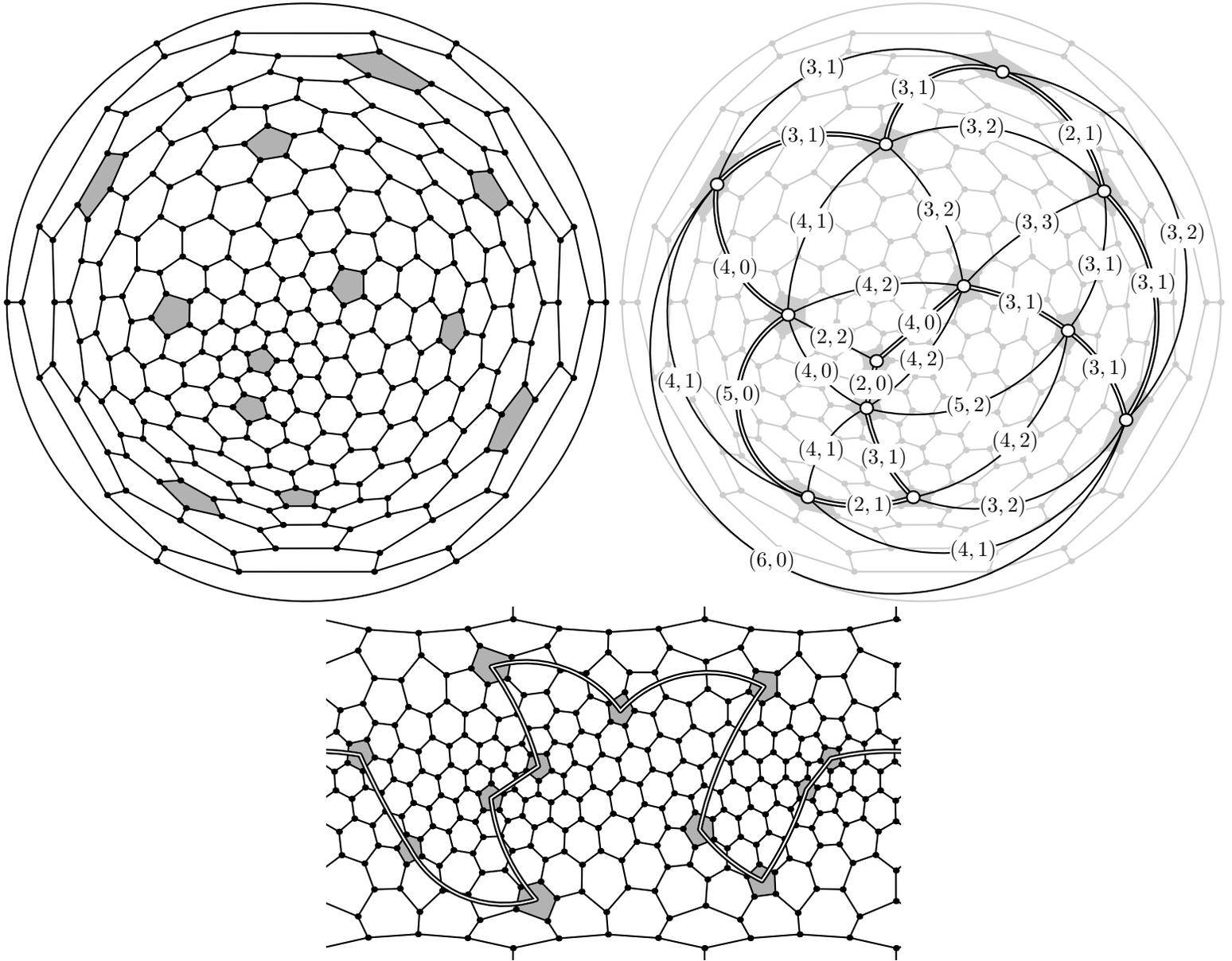

\centerline{
\includegraphics{bigEx.0}
\hfill
\includegraphics{bigEx.1}
}
\centerline{
\includegraphics{bigEx.2}
}
\caption{An example of a Barnette graph on 322 vertices (top left). A triangulation capturing the mutual position of all the small faces with a Hamilton cycle (top right). Another (tubular) drawing of the same graph (bottom); the three edges sticking to the north (to the south) are incident to an omitted vertex at the north (south) pole.}
\label{fig:bigEx}
\end{figure}

In Figure \ref{fig:bigEx}, an example of a Barnette graph on 322 vertices is depicted, along with the corresponding triangulation and a shortest Hamilton cycle in it.

\section{Proof of Theorem \ref{th:find}: Finding a 2-factor}
\label{sec:pro}

In this section we explain the general proceduce in the case when the small faces of $G$ are far from each other.
We will deal with the case when some small faces of $G$ are close to each other in Section \ref{sec:com}.

\subsection{Phase 1: Cut the graph and fix a coloring}

Let $G$ be a Barnette graph, let $T$ be a Hamiltonian triangulation capturing the mutual position of the small faces of $G$, whose existence is given by Lemma \ref{le:triang}. Let $C_T$ be a Hamiltonian cycle in $T$ such that the sum of the lengths of the corresponing Goldberg vectors is minimal. Then there exists a cycle $C^*$ in $G^*$ including all the small vertices of $G^*$ in the same order as the corresponding vertices or $C_T$. 

A cycle in $G^*$ corresponds to an edge-cut in $G$. We cut the graph $G$ along $C^*$. We obtain two graphs, say $G_1$ and $G_2$, containing only hexagons as internal faces, and with semi-edges and partial faces on the boundary.

Both $G_1$ and $G_2$ are subgraphs of the hexagonal grid, hence there is a canonical face coloring using three colors for each of them. We will use colors 1, 2, 3 for one and colors $A$, $B$, $C$ for the other. We color the partial faces in both graphs too.

We choose one color in each graph, say 1 and $A$ (there are 9 color combinations in total), and recolor black all the faces of $G_1$ and $G_2$ colored 1 or $A$; we color white the other faces. (Later we will inspect all the nine colorings.) This gives a black-and-white face-coloring $\phi_i$ inducing a $2$-factor $F_i$ in $G_i$, $i=1,2$.

Observe that for any choice of a color in $G_i$ ($i=1,2$), the edges incident to one face of the other two colors each form a matching $M_i$ such that $G_i^{M_i}= H_i^\parallel$, where $H_i$ is the graph whose vertices are the centers of the faces of the other two colors.

\begin{figure}
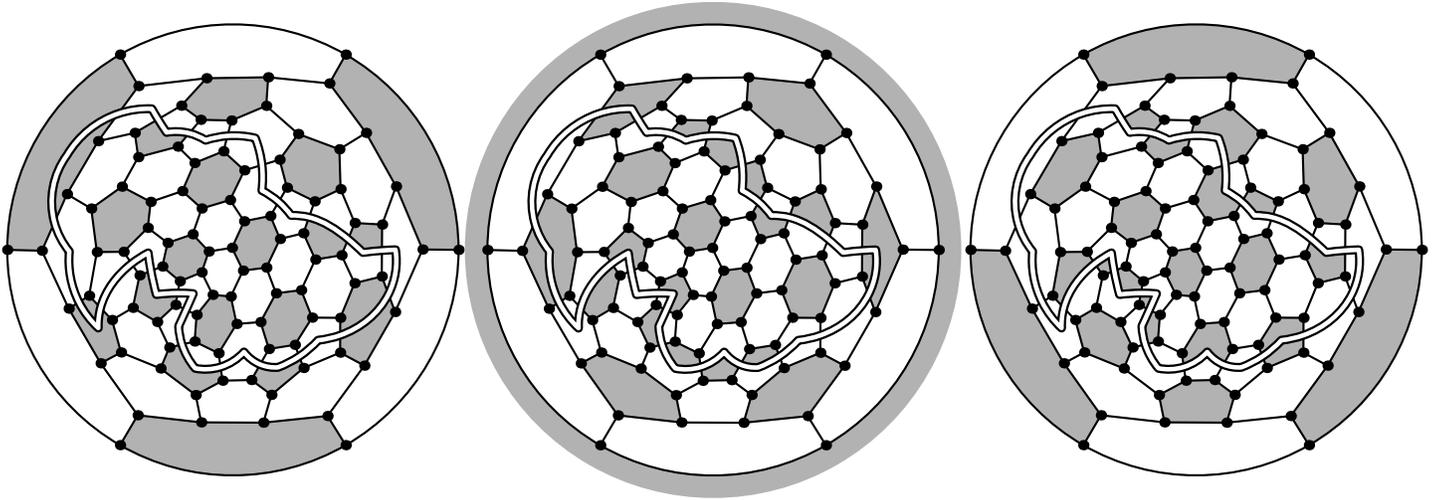

\centerline{
\begin{tabular}{c@{}c@{}c@{}}
\begin{tabular}{@{}c@{}}
\includegraphics{example.2}
\end{tabular}
&
\begin{tabular}{@{}c@{}}
\includegraphics{example.1}
\end{tabular}
&
\begin{tabular}{@{}c@{}}
\includegraphics{example.3}
\end{tabular}
\end{tabular}
}
\caption{Three of the nine black-and-white colorings of the graph in Figure \ref{fig:ex}, (combinations of three different colorings of $G_1$ and three different colorings of $G_2$) corresponding to the given order of pentagons.}
\label{fig:ex2}
\end{figure}

We merge the two black-and-white face-colorings $\phi_1$ and $\phi_2$ of $G_1$ and $G_2$, respectively, into an intermediate black-and-white (multi-)face-coloring $\phi^{(i)}$ of $G$ in the natural way: A face not  corresponding to a vertex of $C^*$ inherits a color from either $G_1$ or $G_2$; A face which is cut by the cycle $C^*$ is divided into two partial faces, one inheriting a color from $G_1$ and the other from $G_2$, see Figure \ref{fig:ex2} for illustration. 

\subsubsection{Active and inactive segments}

The cycle $C^*$ can always be decomposed into a sequence of $\ell\le 12$ subpaths $P^*_1,\dots, P^*_{\ell}$ joining consecutive pairs of small vertices. Let us call these subpaths \emph{segments}. 

We may suppose that a segment only contains hexagons with a non-empty intersection with the straight line joining the end-vertices of the segment.

For each segment $P_i^*$, the two face-colorings of $G_1$ and $G_2$ meet along $P_i^*$, and there is a unique canonical bijection $\varphi_i:\{1,2,3\}\to\{A,B,C\}$ between the two sets of colors. 

If $\varphi_i(1)=A$ then the two black-and-white colorings coincide along $P_i$, we say that the segment $P_i$ is \emph{inactive}; otherwise it is \emph{active}. Out of the nine colorings, each segment is active in precisely six of them.
For example, the segments $p_{6}p_{7}$ and $p_{12}p_1$ are inactive in all the three colorings depicted in Figure \ref{fig:ex2}, the segment $p_4p_5$ is active in all the three colorings, whereas the segment $p_9p_{10}$ is inactive in the first coloring and active in the other two.

\begin{figure}[ht]
\centerline{\includegraphics{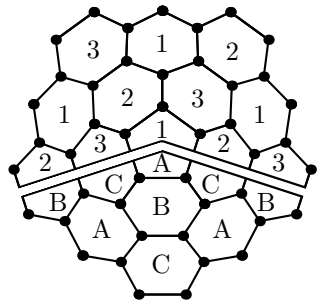}}
\caption{A pentagon always causes a single switch of colors -- the two colors different from its color are switched.}
\label{fig:abc}
\end{figure}

When switching from $P_i^*$ to $P_{i+1}^*$, if the $i$-th small face is a quadrangle, we have $\varphi_i=\varphi_{i+1}$.
If the $i$-th small face is a pentagon, the difference $\varphi_{i+1} \circ \varphi_{i}^{-1}$ is a permutation of the colors $\{A,B,C\}$ such that the color of the pentagon is stable and the two other colors are switched -- a transposition. See Figure \ref{fig:abc} for illustration.

Let $p_i$ be a pentagonal face of $G$ such that the segments $P^*_{i-1}$ and $P^*_i$ meet at $p_i$. Then exactly one of the following happens:

\begin{enumerate}
\item[(i)] if $\varphi_{i-1}(1)=\varphi_i(1)=A$, then both $P^*_{i-1}$ and $P^*_i$ are inactive, $p_i$ generates a switch of $B$ and $C$, thus it is colored $A$ and it is black in both subgraphs;
\item[(ii.a)] if $\varphi_{i-1}(1)=A$ and $\varphi_{i}(1)\ne A$, then $P^*_{i-1}$ is inactive and $P^*_i$ is active, $p_i$ generates a switch of $A$ and $\varphi_{i}(1)$, thus it is colored neither $A$ nor $1$, so it is white in both subgraphs;
\item[(ii.b)] if $\varphi_{i-1}(1)\ne A$ and $\varphi_{i}(1)= A$, then $P^*_{i-1}$ is active and $P^*_i$ is inactive, $p_i$ generates a switch of $A$ and $\varphi_{i-1}(1)$, thus it is colored neither $A$ nor $1$, so it is white in both subgraphs;
\item[(iii.a)] if $\varphi_{i-1}(1)=\varphi_{i}(1)\ne A$, then both $P^*_{i-1}$ and $P^*_i$ are active, $p_i$ generates a switch of $A$ and the third color, thus it is colored $\varphi_{i-1}(1)$, so it is black in $G_1$ and white in $G_2$;
\item[(iii.b)] if $\{\varphi_{i-1}(1),\varphi_{i}(1)\}=\{B,C\}$, then both $P^*_{i-1}$ and $P^*_i$ are active, $p_i$ generates a switch of $B$ and $C$, thus it is colored $A$, so it is white in $G_1$ and black in $G_2$.
\end{enumerate}

In order to transform $\phi^{(i)}$ into a black-and-white face-coloring of $G$ corresponding to a good $2$-factor of $G$, we reroute slightly the cut $C^*$ in a way described in the following subsection.

\subsection{Phase 2: Approximate the cut by $\Gamma$-paths}

Let $P^*_i$ be an active segment, let $\varphi_i(1)=B$. Suppose without loss of generality that $\varphi_i^{-1}(A)=2$.
Then all the faces of $P_i$ colored $A$ (and 2) or $1$ (and $B$) are partially black and partially white; both parts of each face of $P_i^*$ colored $C$ and 3 are white. 

We approximate the dual path $P^*_i$ by a sequence $Q_i$ of faces colored $C$ and/or 3, each consecutive pair of faces in a mutual position $(1,1)$.

Let $f$ be a white ($C$- and $3$-colored) hexagonal face of $Q_i$. Then among its neighbors, there is a cyclic sub-sequence of $A$- and $B$-faces colored alternatively black and white, and another cyclic sub-sequence of $1$- and $2$-faces colored alternatively black and white, with the coloring being the opposite of the first one. Therefore, there are exactly two pairs (not necessarily disjoint) of adjacent faces of the same color: each pair is either a black $A$-face adjacent to a black $1$-face, or a white $B$-face adjacent to a white $2$-face. Therefore, $f$ is a white non-resonant hexagon, corresponding to a vertex of degree 2 in the future auxiliary graph $\Gamma$ being constructed -- we will call it a $\Gamma$-face.

Let $f$ and $f^\prime$ be two consecutive $\Gamma$-faces. If the two faces adjacent both to $f$ and $f^\prime$ are black, then the two cycles of the $2$-factors in $G_1$ and $G_2$ are merged. If the two faces adjacent both to $f$ and $f^\prime$ are white, then a new 2-cycle of the $2^*$-factor is created. In the first case, the $\Gamma$-edge $ff^\prime$ is black, in the second case it is a white one.

Two consecutive $\Gamma$-edges of $Q_i$ of the same color always form a $180^\circ$ angle, otherwise it could be possible to simplify $Q_i$ by removing a face from $Q_i$. Similarly, two consecutive edges of $Q_i$ of different colors always form an angle of  $\pm 120^\circ$.

\begin{figure}[ht]
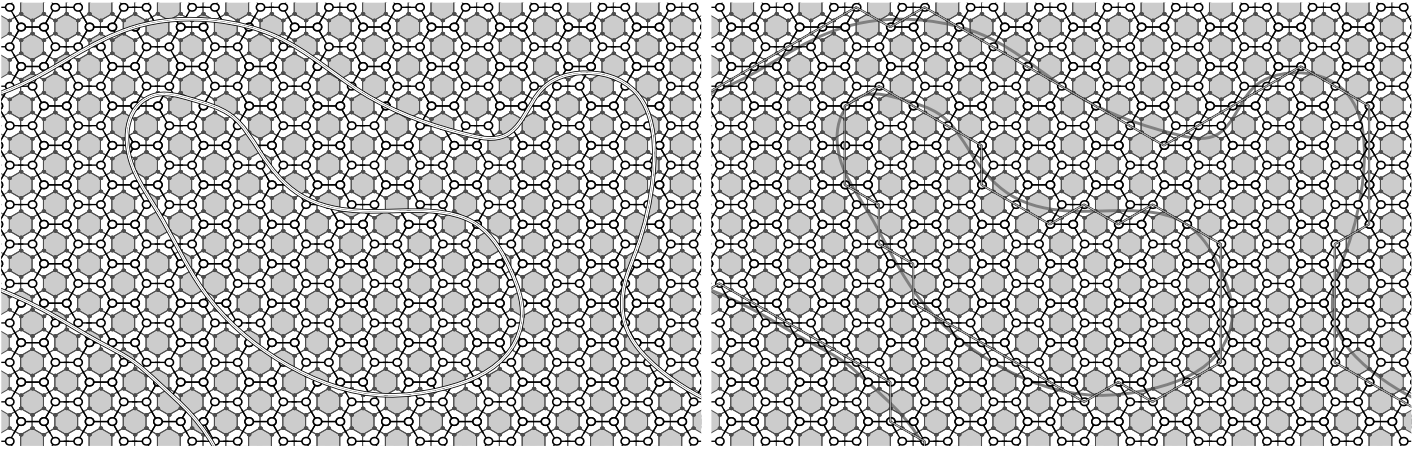

\centerline{
\includegraphics[scale=.35]{grid.0}
\hfil
\includegraphics[scale=.35]{grid.1}
}
\caption{Several hexagonal patterns meeting along some cut curves (left). As the cutting lines are approximated by $\Gamma$-paths, a good $2^*$-factor is created (right).
}
\label{fig:grid1}
\end{figure}

%
The resulting structure of $\phi^{(i)}$ along $P^*_i$ is the following: All vertices are covered by cycles of length 6 (single faces), 10 (two adjacent black hexagons, both incident to a black $\Gamma$-edge), or 2 (white $\Gamma$-edges). A $\Gamma$-path $Q_i$ separates the two subgraphs of regular coloring. See Figure \ref{fig:grid1} for illustration.

The first (the last) $\Gamma$-face of $Q_i$ is the pentagon $p_i$ ($p_{i+1}$) if and only if the segment $P^*_{i-1}$ ($P^*_{i+1}$) is inactive; otherwise the first (the last) $\Gamma$-face of $Q_i$ is a hexagon adjacent to $p_i$ ($p_{i+1}$) and it is the last (the first) $\Gamma$-face of $Q_{i-1}$ ($Q_{i+1}$, respectively). 

White non-resonant hexagons where two consecutive sequences $Q_{i-1}$ and $Q_i$ meet are the only occasion where two $\Gamma$-edges of the same color might form a $60^\circ$ angle -- if only they are both incident to the same pentagon.

Let us explicit the structure of $H^{(i)}=H_1\cup H_2$ and of $\Gamma$ now:
Vertices of $H^{(i)}$ are all the vertices corresponding to faces of $G$ white in $G_1$ or in  $G_2$; 
each vertex of $\Gamma$ where two black edges meet corresponds to a 2-vertex in $H^{(i)}$ (the corresponding face of $G$ is a non-resonant white hexagon adjacent to four black faces belonging to two different components of the $2^*$-factor); each vertex of $\Gamma$ where two white edges meet corresponds to a 4-vertex in $H^{(i)}$ (the corresponding face of $G$ is incident to four different compents of the $2^*$-factor, including two 2-cycles); each vertex of $\Gamma$ where a black and a white edge meet at a $120^\circ$ angle corresponds to a 3-vertex in $H^{(i)}$ (the corresponding face of $G$ being incident to three different components of the $2^*$-factor: a 2-cycle, a 6-cycle and a 10-cycle). 

If there are $q_5$ white pentagons, then $\Gamma$ is composed of $q_5/2$ paths. A white quadrangle is either an isolated vertex of $\Gamma$ (if both incident segments are inactive) or it is an internal vertex of a path (otherwise).

\subsection{Phase 3: Change the parity of the $2^*$-factor}

It follows from Lemma \ref{l:paths} that whenever we want to transform an even $2^*$-factor into an odd one, it suffices either to increase or decrease the number of black quadrangles by 1, or to increase or decrease the number of black pentagons by 2. In other words, it suffices either to change the number of isolated vertices in $\Gamma$ by 1 or change the number of $\Gamma$-paths by 1.

\subsubsection{Changing the parity using a quandrangle} 

Let $q$ be a quadrangular face of $G$. For three of the nine colorings of $G_1$ and $G_2$, both segments incident to $q$ are inactive; moreover, for two out of the three $q$ is a white face. In Phase 1, we choose one of these two.

\begin{figure}[ht]
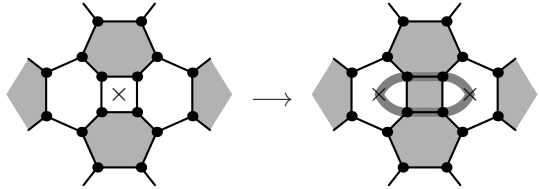

\centerline{
\begin{tabular}{@{}c@{}c@{}c@{}}
\begin{tabular}{c}
\includegraphics{quad.5}
\end{tabular}
&
$\longrightarrow$
&
\begin{tabular}{c}
\includegraphics{quad.6}
\end{tabular}
\end{tabular}
}
\caption{We can use a white quadrangle to change the parity of a $2^*$-factor. The times sign marks non-resonant faces -- vertices of $\Gamma$; edges of $\Gamma$ are drawn using a thick grey line.}
\label{fig:quad}
\end{figure}

If the good $2^*$-factor obtained in Phase 1 is even, it can be transformed into an odd one by recoloring $q$ black. This way an isolated vertex of $\Gamma$ is transformed into a cycle of length 2, see Figure \ref{fig:quad} for illustration.

\subsubsection{Changing the parity using two pentagons}

From this point on we may assume that $G$ has no quadrangular faces -- it is a (fullerene) graph having 12 pentagonal faces. 

Suppose first that some pair of consecutive pentagons $p_i$ and $p_{i+1}$ (consecutive along the cut $C$) are in the mutual position $(c_1,c_2)$, $c_1\ge c_2\ge 0$, with $3\mid (c_1-c_2)$. Then in the coloring of $G_1$ with colors 1, 2, 3 (and of $G_2$ with $A$, $B$, $C$) the partial faces corresponding to the pentagons $p_i$ and $p_{i+1}$ have the same color. Therefore, for two of the nine colorings the segment $P_i$ joining $p_i$ and $p_{i+1}$ is active whereas the neighboring segments $P_{i-1}$ and $P_{i+1}$ are inactive.

For both such colorings, after Phase 2 there is a $\Gamma$-path with endvertices at $p_i$ and $p_{i+1}$, and the vertex set of this path can be chosen to be the same in both colorings. If this is the case, then each $\Gamma$-edge white in one coloring is black in the other and vice versa. Among the two colorings, we may fix the one where the number of white $\Gamma$-edges is maximised.

\begin{figure}[ht]
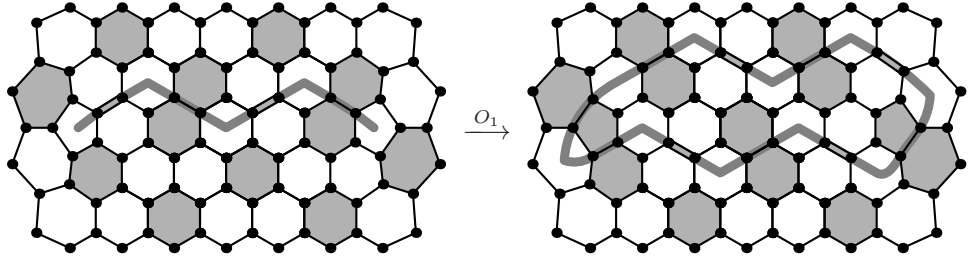

\centerline{
\begin{tabular}{@{}c@{}c@{}c@{}}
\begin{tabular}{c}
\includegraphics{pathtocycle.1}
\end{tabular}
&$\xlongrightarrow{O_1}$&
\begin{tabular}{c}
\includegraphics{pathtocycle.2}
\end{tabular}
\end{tabular}
}
\caption{Operation $O_1$: The parity of a $2^*$-factor can be changed by modifying a $\Gamma$-path joining two consecutive pentagons into a $\Gamma$-cycle.}
\label{fig:pc}
\end{figure}

We transform the $\Gamma$-path into a $\Gamma$-cycle, increasing the number of black pentagons by 2, in the following way:
For each black $\Gamma$-edge, we recolor both black hexagons forming a black 10-cycle white; then we recolor all faces corresponding to the vertices of the $\Gamma$-path black, including the first and the last one ($p_i$ and $p_{i+1}$). We will denote this operation $O_1$. See Figure \ref{fig:pc} for illustration.

\begin{figure}[ht]
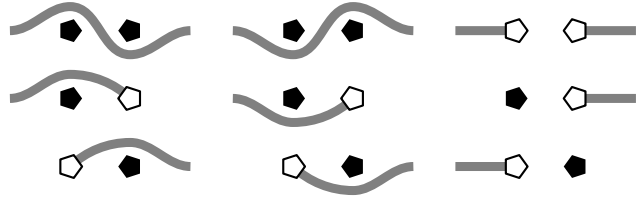

\centerline{
\begin{tabular}{ccc}
\includegraphics[scale=1]{schema.10}
&\includegraphics[scale=1]{schema.11}
&\includegraphics[scale=1]{schema.18}
\\\includegraphics[scale=1]{schema.13}
&\includegraphics[scale=1]{schema.12}
&\includegraphics[scale=1]{schema.14}
\\\includegraphics[scale=1]{schema.16}
&\includegraphics[scale=1]{schema.15}
&\includegraphics[scale=1]{schema.17}
\end{tabular}
}
\caption{A schematic drawing of the position of the $\Gamma$-paths in the neighborhood of two consecutive pentagons of different colors.}
\label{fig:9cases}
\end{figure}

From this point on we may assume that there is no pair of consecutive pentagons with the same color in $G_1$ (or in $G_2$). Then for every pair of consecutive pentagons the nine colorings look like depicted in Figure \ref{fig:9cases}.

Let $\phi_i^j$ be the angle between the two segments meeting at pentagon $p_i$ in $G_j$, $j=1,2$. Clearly, $\phi_i^1+\phi_i^2=300^\circ$. When following the segments composing the cut in an ascending order, say $G_1$ is to the left and $G_2$ to the right. If $\phi_i^1 > 150^\circ > \phi_i^2$, then there is a right turn at $p_i$ when switching from $P_{i-1}$ to $P_i$. If $\phi_i^1 < 150^\circ < \phi_i^2$, then there is a left turn at $p_i$ when switching from $P_{i-1}$ to $P_i$. The value $\phi_i^1=\phi_i^2=150^\circ$ means that the segment $P_i$ continues in the same direction as $P_{i-1}$. 

Let $\phi_i=\phi_i^1-\phi_i^2$ for $i=1,\dots,12$. It is easy to see that $\sum_{i=1}^{12}\phi_i=0$, since $\sum_{i=1}^{12}\phi_i^1 = \sum_{i=1}^{12}\phi_i^2 = 1800^\circ$.  Therefore, there exist $i$ such that $\phi_i\cdot \phi_{i+1}\le 0$ (indices modulo 12). We fix $i$ such that $\phi_i\cdot\phi_{i+1}\le 0$ and the difference $|\phi_i-\phi_{i+1}|$ is as big as possible. 

Without loss of generality we may assume that $\phi_i\ge 0$ and $\phi_{i+1}\le 0$. In other words, there is a right turn at $p_i$ followed by a left turn at $p_{i+1}$. There are two colorings in which the segments $P_{i-1}$, $P_i$, and $P_{i+1}$ are active; among them we choose the one where $p_i$ is black in $G_2$ and $p_{i+1}$ is black in $G_1$.

\begin{figure}[ht]
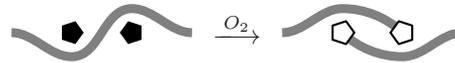

\centerline{
\begin{tabular}{@{}c@{}c@{}c@{}}
\begin{tabular}{c}
\includegraphics{schema.11}
\end{tabular}
&$\xlongrightarrow{O_2}$&
\begin{tabular}{c}
\includegraphics{schema.19}
\end{tabular}
\end{tabular}
}
\caption{A schematic drawing of the operation $O_2$.}
\label{fig:o2}
\end{figure}

\begin{figure}[ht]
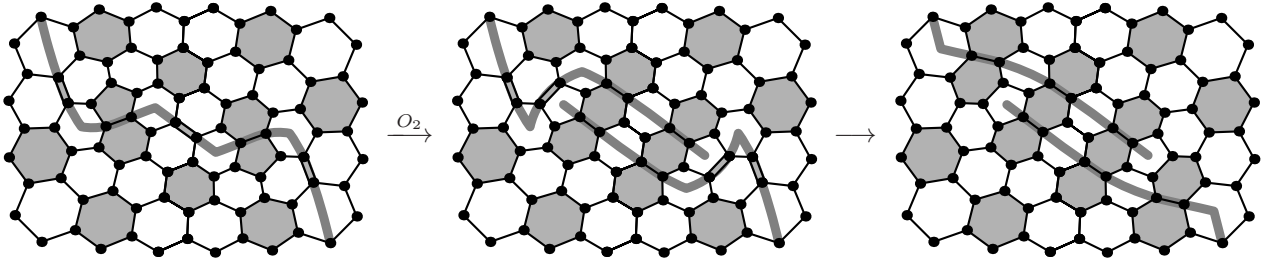

\centerline{
\begin{tabular}{@{}c@{}c@{}c@{}c@{}c@{}}
\begin{tabular}{c}
\includegraphics{transform.3}
\end{tabular}
&$\xlongrightarrow{O_2}$&
\begin{tabular}{c}
\includegraphics{transform.4}
\end{tabular}
&$\longrightarrow$&
\begin{tabular}{c}
\includegraphics{transform.5}
\end{tabular}
\end{tabular}
}
\caption{Operation $O_2$: The parity of a $2^*$-factor can be changed by transforming a $\Gamma$-path passing by two consecutive pentagons into two different $\Gamma$-paths.}
\label{fig:tr1}
\end{figure}

\begin{figure}[pht]
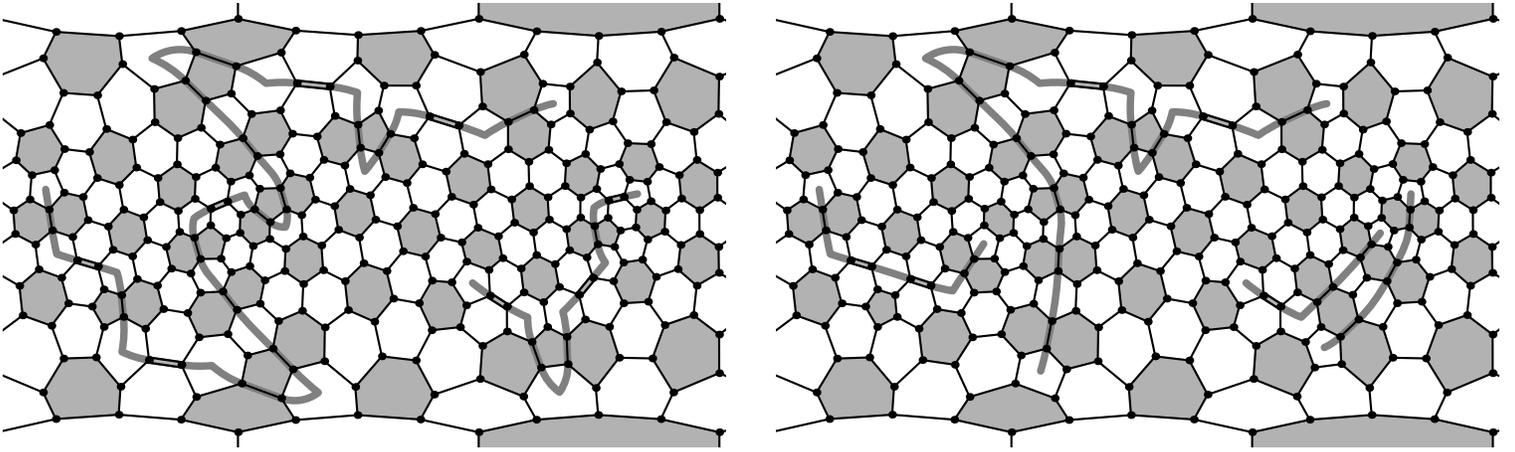

\centerline{
\includegraphics{bigEx.3}
\hskip 5mm
\includegraphics{bigEx.30}
}
\caption{The good $2^*$-factor induced by one of the nine possible black-and-white face-colorings of the graph in Figure \ref{fig:bigEx} (left). It is already an odd $2^*$-factor; there are several pairs of pentagons for which the operation $O_2$ is admissible. Another good $2^*$-factor of the same graph obtained by two applications of $O_2$ (right).}
\label{fig:bigEx2}
\end{figure}

We can now change the parity of the $2^*$-factor (if needed) by decreasing the number of black pentagons in the following way: For each black $\Gamma$-edge of $P_i$, we recolor both black hexagons forming a black 10-cycle white; then we recolor all faces corresponding to the vertices of the $\Gamma$-subpath $Q_i$ black, including the first and the last one (those adjacent to $p_i$ and $p_{i+1}$, respectively); we recolor $p_i$ and $p_{i+1}$ white. As the last step, we simplify unnecessary $60^\circ$ turns. We will denote this operation $O_2$. See Figures \ref{fig:o2},
\ref{fig:tr1} and \ref{fig:bigEx2} for illustration.

\subsection{Phase 4: Transform a good odd $2^*$-factor into a simple $2$-factor}

It suffices now, as the last phase, to transform a good odd $2^*$-factor into a simple (odd) $2$-factor. We do it in the following way:

In a good $2^*$-factor, each 2-cycle corresponds to a white $\Gamma$-edge $ff^\prime$, incident to two white resonant hexagons $h_1$ and $h_2$ (one in each of $G_1$ and $G_2$). We can choose either $h_1$ or $h_2$, say $h_i$, and recolor it black: By doing this, the 2-cycle is merged with two other cycles in $G_i$; the other face $f_0$ incident to both cycles being merged loses its resonantness, it becomes another $\Gamma$-face inserted to the $\Gamma$-path between $f$ and $f^\prime$, joint now to $f$ and $f^\prime$ by two black $\Gamma$-edges forming a $60^\circ$ angle and replacing the original white $\Gamma$-edge. In $H^{(i)}$, a vertex of degree 3 is removed, and thus the degree of three other vertices is decreased by 1: one of them corresponds to $f_0$, the other two correspond to $f$ and $f^\prime$. 

Observe that this operation decreases the number of components of the factor by 2, therefore, starting with an odd factor we can only obtain odd factors.

We make a decision for all white $\Gamma$-edges sequentially according to their order along $Q_i$, according to the following rules: If a white $\Gamma$-edge $e_j$ forms a $180^\circ$ angle with $e_{j-1}$ (which has to have been white in this case) and that we have decided to recolor black a hexagon in $G_i$, $i=1,2$, incident to $e_{j-1}$, then we decide to recolor black a hexagon in $G_{3-i}$ incident to $e_j$. If a white $\Gamma$-edge $e_j$ forms a $120^\circ$ angle with a black $e_{j-1}$, we decide to recolor black a hexagon incident to $e_j$ in such a way that one of the new black $\Gamma$-edges forms a $180^\circ$ angle with $e_{j-1}$.

The resulting structure in $G$ is the following: All the $\Gamma$-paths and $\Gamma$-cycles are
formed of black $\Gamma$-edges only. Each vertex of $\Gamma$ of degree 1 or 2 corresponds to a 2-vertex in $H^{(i)}$.

Finally, to obtain $H$, we suppress all the 2-vertices in $H^{(i)}$; for each $\Gamma$-edge we merge the incident partial faces of $H^{(i)}$.

To describe the structure of $H$, we introduce the following notation:
A vertex of $\Gamma$ is called \emph{direct} if it corresponds to a pentagon or if the two incident (black) $\Gamma$-edges form a $180^\circ$ degree; otherwise it is called \emph{sharp}. 

We claim that there cannot be three consecutive sharp $\Gamma$-vertices along any $Q_i$: Suppose some $Q_i$ contains a subpath $f_0f_1f_2f_3f_4$ with all of $f_1$, $f_2$, and $f_3$ sharp and $f_0$ direct. If $f_1f_3$ had been a white $\Gamma$-edge after the Phase 2, we would not have decided to choose $f_2$. Therefore, $f_2$ was a $\Gamma$-vertex already after Phase 2, which means that $f_0f_2$ was a white $\Gamma$-edge after Phase 2. If $f_2f_4$ was also a white $\Gamma$-edge after Phase 2, we would have decided one of them in the other way. Therefore, $f_3$ was a $\Gamma$-vertex already after Phase 2, but not $f_4$, which means that $f_4$ is sharp. As $f_1$ must have been chosen because of the other $\Gamma$-edge incident to $f_0$, $f_4$ should never have been chosen, a contradiction.

A (black) $\Gamma$-edge joining two direct $\Gamma$-vertices $f$ and $f^\prime$ completes the boundary of two partial faces in $H_1$ and $H_2$, each having three incident 3-vertices. After the suppression of $2$-vertices in $H^{(i)}$, in $H$ these two partial faces are merged into a hexagon.

The $60^\circ$ angle at a sharp $\Gamma$-vertex $f$ contains a partial face of $H^{(i)}$ having one 3-vertex, which is to be merged with (at least) two other partial faces. 

If both $\Gamma$-vertices adjacent to $f$ in $\Gamma$ are direct, then a face of size 7 is created in $H$ by merging two partial faces each having three incident 3-vertices in $H_i$ with a partial face having one incident 3-vertex in $H_{3-i}$.
On the other hand, opposite to this one, there is a face of $H^{(i)}$ whose size is decreased by 1 by the suppresion of the 2-vertex $f$ -- a pentagonal face is created in $H$.

\begin{figure}[pht]
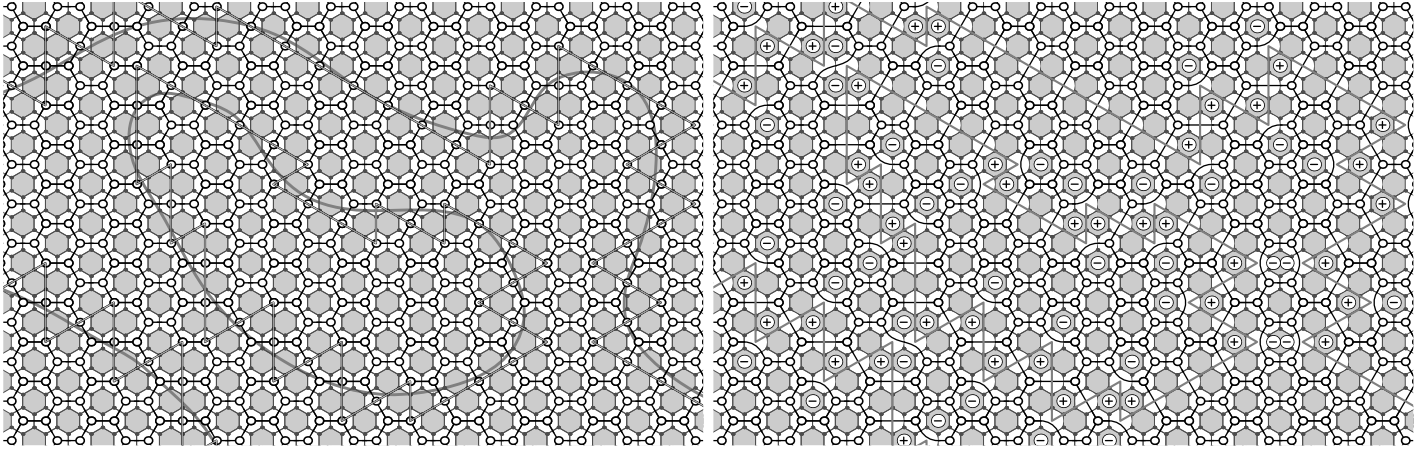

\centerline{
\includegraphics[scale=.35]{grid.2}
\hfil
\includegraphics[scale=.35]{grid.3}
}
\caption{The intermediate structure after eliminating the white $\Gamma$-edges (left) and the final simple $2$-factor, with the face size changes in $H$ (with respect to the initial size of 6) marked with plus and minus signs.
}
\label{fig:grid2}
\end{figure}

\begin{figure}[pht]
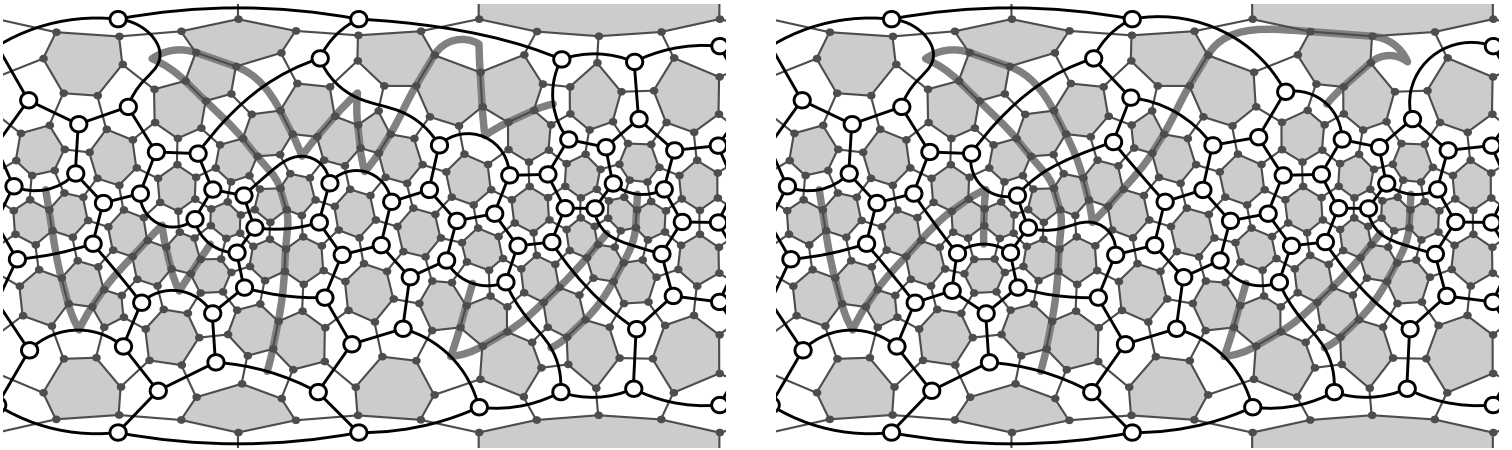

\centerline{
\includegraphics{bigEx.31}
\hskip 5mm
\includegraphics{bigEx.32}
}
\caption{The simple $2$-factor obtained from the good $2^*$-factor in Figure \ref{fig:bigEx2}, depicted together with the residual graph $H$ and the auxilliary graph $\Gamma$ (left). Another simple $2$-factor obtained from the previous one by "flipping out" unnecessary zig-zags of sharp $\Gamma$-vertices (right). Observe that for the latter, the residual graph $H$ has faces of size 5, 6, and 7 only and it is cyclically 5-edge-connected.}
\label{fig:bigEx3}
\end{figure}

\begin{figure}[pht]
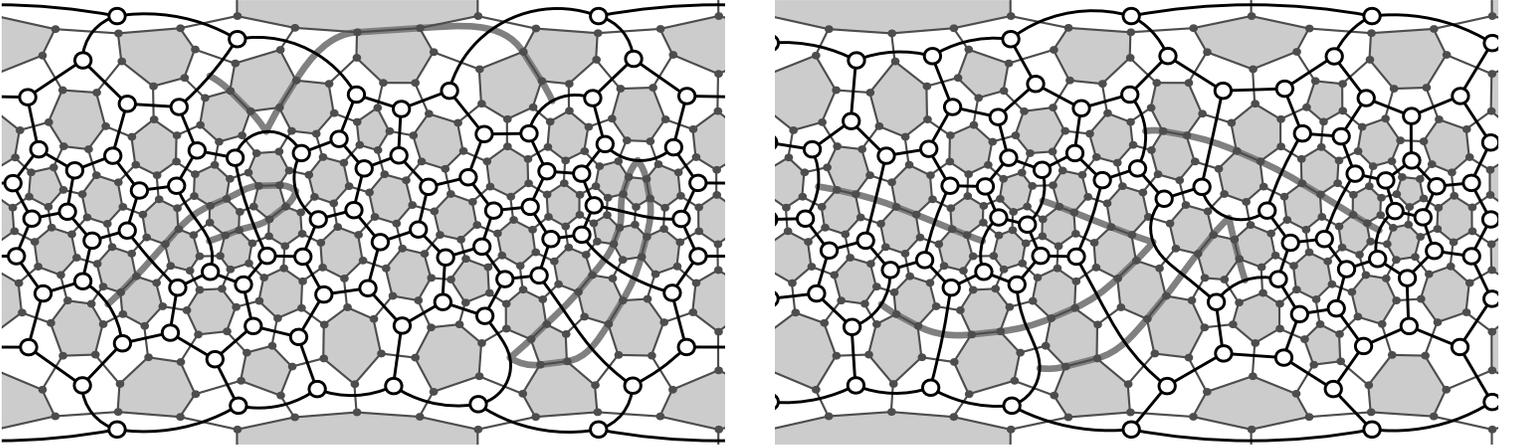

\centerline{
\includegraphics{bigEx.7}
\hskip 5mm
\includegraphics{bigEx.8}
}
\caption{Two different odd simple $2$-factors of the graph in Figure \ref{fig:bigEx} with the largest number of vertices of the residual graph $H$ (82) we were able to find.}
\label{fig:bigEx4}
\end{figure}

If one of the vertices adjacent to a sharp vertex in $\Gamma$ is a sharp one, they are transformed  into a face of size 8 and two pentagons in $H$. See Figures \ref{fig:grid2}, \ref{fig:bigEx3}, and  \ref{fig:bigEx4} for illustration.

\section{Checking the correctness of the algorithm in the neighborhood of small faces close to each other: the computer-assisted part}
\label{sec:com}

Let $G$ be a Barnette graph. Let $S(G)$ be the set of the \emph{small} faces (faces of size 4 or 5) of $G$. It is straightforward to derive from the Euler's formula that $2f_4+f_5=12$, where $f_4$ and $f_5$ are the numbers of quadrangles and pentagons in $G$, respectively. 

\subsection{Patches}

A \emph{patch} is a 2-connected subcubic plane graph $P$, having at most one face of size different from 4, 5 and 6, and such that all vertices of $P$ of degree 2 are incident to this special face, often referred to as the outer face of the patch; moreover, $P$ contains no pair of adjacent 4-faces. When a patch is depicted, there are additional pending half-edges at vertices of degree 2 towards the outer face.

The \emph{curvature} of a patch $P$, denoted by $\mu(P)$, is equal to $2f_4(P)+f_5(P)$, where $f_4(P)$ and $f_5(P)$ are the numbers of quadrangles and pentagons in $P$ (distinct from the outer face of $P$), respectively. 

We denote $\partial(P)$ the \emph{boundary} of a patch $P$ -- the facial cycle of the outer face of $P$; we denote $\delta(P)$ the \emph{perimeter} of a patch $P$, the number of 2-vertices in $P$.

The \emph{boundary vector} $\sigma(P)$ of a patch $P$ is a cyclic sequence of distances between consecutive 2-vertices on the boundary cycle of $P$. The length of $\sigma(P)$ is equal to $\delta(P)$ and its sum is equal to the length of $\partial(P)$. When expliciting elements of a cyclic sequence $\sigma$, we write $x^k$ as a shortcut for $k$ consecutive occurences of a value $x$ in $\sigma$.

 Each vertex of $\partial(P)$ is either a 2-vertex or a 3-vertex in $P$. The proportion of 2-vertices along $\partial(P)$ is determined by the curvature of $P$, as is stated explicitely in the following lemma, which is a generalisation of an observation from \cite{KS} and can be derived directly from Euler's formula by the same double-counting arguments.
 
\begin{lemma}
Let $P$ be a patch of curvature $\mu$. Then 
$$
2\delta(P)-|\partial(P)| = 6-\mu.
$$
\label{lemma:basic}
\end{lemma}
Observe that for patches of curvature (greater than, less than) six, the average value of $\sigma(P)$ is (greater than, less than, respectively) two.

A patch $P$ of curvature $\mu\le 4$ ($\mu\ge 8$) is called \emph{convex} if its boundary vector $\sigma(P)$ only contains '1's and '2's ('2's and '3's, respectively). A patch $P$ with $\mu=5$ ($\mu=7$) is called \emph{convex} if $\sigma(P)$ does not contain $32^j3$ (does not contain $12^j1$, respectively). 
A patch $P$ of curvature 6 is called \emph{convex} if $\sigma(P)$ contains at most one subsequence $32^j3$; if this is the case, $j>\delta(P)/2$. For instance, all the caps of nanotubes in Figures \ref{fig:33tube} and \ref{fig:caps} are convex patches of curvature 6. 

Note that, according to Lemma \ref{lemma:basic}, the boundary vector of a convex patch $P$ of curvature $\mu\le 4$ has the form $(12^{k_1}12^{k_2}\dots12^{k_{t}})$
where $t=6-\mu$, $k_1,k_2,\dots,k_{t}\in \mathbb{N}_0$, and $k_1+k_2+\dots+k_{t}=p-t$.

\begin{figure}[th]
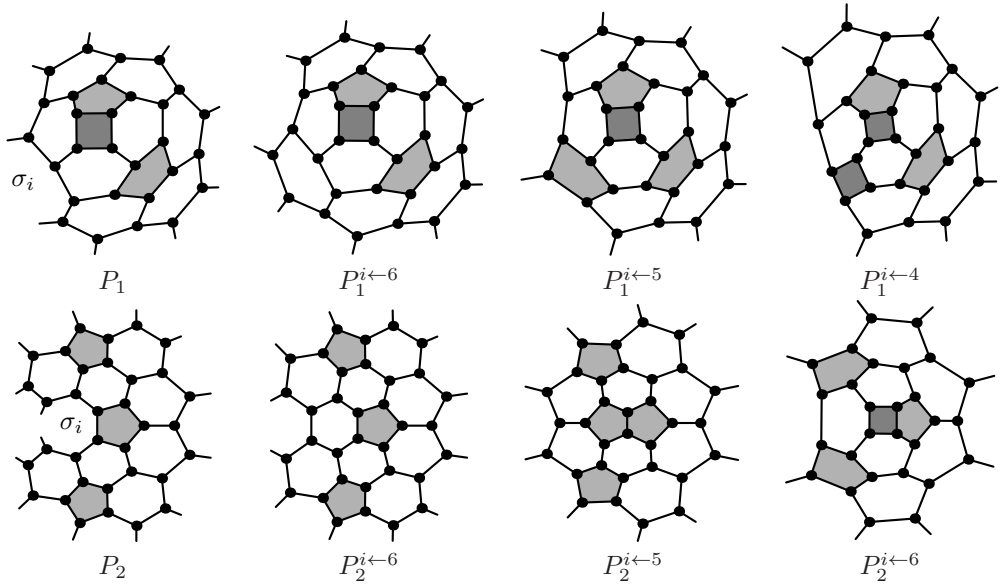

\centerline{
\begin{tabular}{cccc}
\includegraphics{grow.0}&
\includegraphics{grow.6}&
\includegraphics{grow.5}&
\includegraphics{grow.4}\\
$P_1$ & $P_1^{i\gets 6}$ & $P_1^{i\gets 5}$ & $P_1^{i\gets 4}$ \\
\includegraphics{grow.1}&
\includegraphics{grow.16}&
\includegraphics{grow.15}&
\includegraphics{grow.14}\\
$P_2$ & $P_2^{i\gets 6}$ & $P_2^{i\gets 5}$ & $P_2^{i\gets 6}$
\end{tabular}
}
\caption{Two different examples of three different patches obtained from a given patch (on the left) by inserting a new face at the element $\sigma_i$ of its boundary vector. 
}
\label{fig:grow}
\end{figure}

We denote $P^{i\gets j}$ a patch obtained from $P$ by adding a face of size $j$ to $P$ along the path corresponding to the $i$-th element of $\sigma(P)$, if such a patch exists, see Figure \ref{fig:grow} for illustration. 

It may happen that while adding a new face to a patch, we have to identify some elements (vertices/edges/faces) of the patch, as in the second row of Figure \ref{fig:grow}. It may even happen that adding a new face of some desired size to a specific place of a patch is not possible, since the faces to be identified are not of the same size.

\subsubsection{Patches in Barnette graphs}

Let $G$ be a Barnette graph. We say that a patch $P$ is contained in $G$ if there is a graph homomorphism $\varphi: P \to G$ such that all faces of the patch (except for the outer face) are also faces of $G$. 
We say that a patch $P$ is \emph{realizable} if it is contained in some Barnette graph.

Observe that a patch $P$ of perimeter 0 is contained in a Barnette graph $G$ if and only if $P=G$ and the outer face of $P$ is a face of $G$. Similarly, a patch $P$ of perimeter 2 is contained in a Barnette graph $G$ if and only if $P=G\setminus e$ for some edge $e$ of $G$ and the outer face of $P$ is the union of the two faces incident to $e$ in $G$.
Finally, since Barnette graphs are cyclically 4-edge-connected, a patch $P$ if perimeter 3 is contained in a Barnette graph $G$ if and only of $P=G\setminus v$ for some vertex $v$ of $G$ and the outer face of $P$ is the union of the three faces incident to $v$ in $G$.
On the other hand, no patch of perimeter 1 can be realizable, since it would correspond to a cut-edge in a Barnette graph.

Some (but not all) realizable patches can be obtained in the following way: For any induced cycle $C$ of a Barnette graph $G$, there are two distinct (but not disjoint) patches $P$ and $\bar{P}$ contained in $G$ such that $\partial(P)=\partial(\bar{P})=C$. 
It is easy to see that we have $\mu(P)+\mu(\bar{P})=12$ and that $\delta(P)$ is equal to the number of edges of the cut separating $P$ from $G\setminus P$. 

Moreover, as each vertex of $C$ is either a 2-vertex in $P$ or a 2-vertex in $\bar{P}$, $\delta(P)+\delta(\bar{P})$ is equal to the length of $\partial(P)$. 
 
As a direct consequence of Lemma \ref{lemma:basic} we obtain the following observation.	
\begin{lemma}
Let $C$ be an induced cycle in a Barnette graph and let $P$ and $\bar{P}$ be the two corresponding patches. Then 
$$
\delta(P)-\delta(\bar{P}) = 6-\mu(P) = 6+\mu(\bar{P}).
$$
\label{lemma:basic2}
\end{lemma}
However, there are patches contained in Barnette graphes which cannot be obtained this way: it is not always true that the facial cycle of the outer face of a patch corresponds to an induced cycle of the host Barnette graph -- a patch can even be self-overlapping. 

\begin{lemma}
Let $P$ be a realizable patch of perimeter at least $2$. For every element $\sigma_i$ of its boundary vector there exists $j\in \{4,5,6\}$ such that $P^{i\gets j}$ is also a realizable patch, moreover, it is contained (at least) in the same Barnette graph as $P$.
\label{le:real}
\end{lemma}

Proof. Let $P$ be contained in a Barnette graph $G$. Each element of $\sigma$ is a path contained in a facial cycle of some face of $G$ of a certain size $j\in \{4,5,6\}$.  Therefore, the face added to the patch corresponds to a face of $G$.
\ep

\subsubsection{Primitive patches}

\begin{figure}[ht]
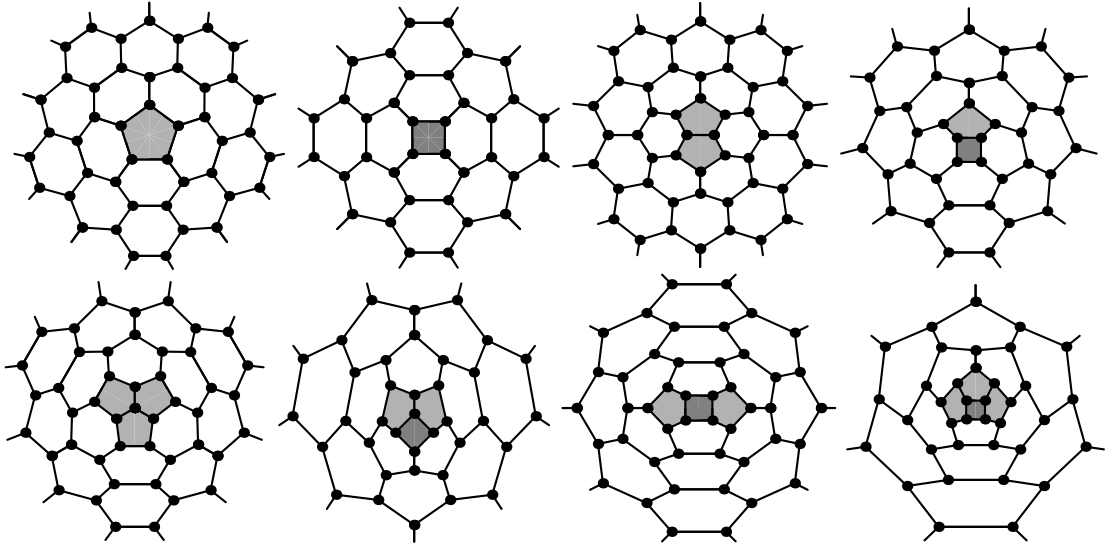

\centerline{
\includegraphics[scale=1.0]{smp.1}
\includegraphics[scale=1.0]{smp.2}
\includegraphics[scale=1.0]{smp.3}
\includegraphics[scale=1.0]{smp.4}
}
\centerline{
\includegraphics[scale=1.0]{smp.5}
\includegraphics[scale=1.0]{smp.6}
\includegraphics[scale=1.0]{smp.7}
\includegraphics[scale=1.0]{smp.8}
}
\caption{Arrangements of small faces in primitive patches for different values of curvature $1\le \mu \le 5$.
}
\label{fig:smallpatches1}
\end{figure}

A convex patch of curvature $\mu\le 5$ is \emph{primitive} if the arrangement of its small faces is the same as in one of the patches depicted in Figure \ref{fig:smallpatches1} or it has no small faces at all (for $\mu=0$). 

Observe that each convex patch with at most one small face is primitive.

\begin{lemma}
Let $P$ be a convex patch of curvature $\mu\le 5$ which is not primitive. Then there exists another patch $P^\prime$ with the same curvature and the same boundary vector as $P$ on a bigger number of vertices.
\label{le:primitive}
\end{lemma}

Proof.
Suppose that there exists a convex patch of curvature $\mu\le 5$ which is not primitive, and all the convex patches of given curvature and boundary vector have at most as many vertices as $P$. 

If $P$ has at most one small face, then it is primitive by definition, a contradiction. Therefore, we may assume that $P$ has at least two small faces. 

If all the small faces of $P$ are pairwise adjacent to each other, then $P$ has at most three small faces, moreover, if it has three small faces, at most one of them is a quadrangle. In all the cases the patch is primitive, a contradiction. 

\begin{figure}[ht]
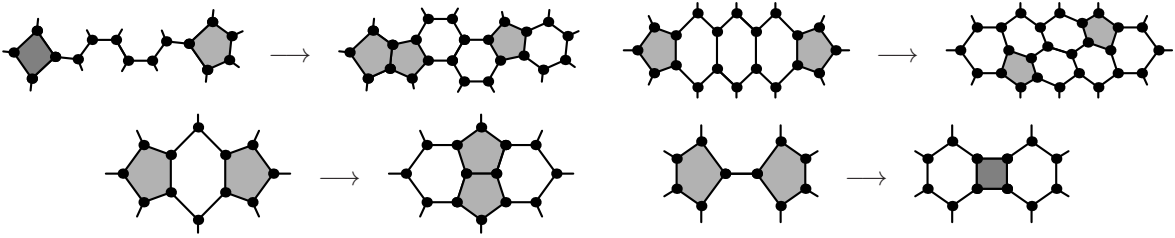

\centerline{
\begin{tabular}{c}
\includegraphics{oper.2}
\end{tabular}
$\longrightarrow$
\begin{tabular}{c}
\includegraphics{oper.1}
\end{tabular}
\begin{tabular}{c}
\includegraphics{oper.8}
\end{tabular}
$\longrightarrow$
\begin{tabular}{c}
\includegraphics{oper.7}
\end{tabular}
}
\centerline{
\begin{tabular}{c}
\includegraphics{oper.4}
\end{tabular}
$\longrightarrow$
\begin{tabular}{c}
\includegraphics{oper.3}
\end{tabular}
\hfil
\hfil
\begin{tabular}{c}
\includegraphics{oper.6}
\end{tabular}
$\longrightarrow$
\begin{tabular}{c}
\includegraphics{oper.5}
\end{tabular}
}
\caption{If a patch contains at least two non-adjacent small faces, it can be transformed to another one with more vertices, unless the two small faces are in position (1,1) and the edge connecting them is incident to a quadrangle: Two generic cases (top) and two special cases (bottom). The size of the two small faces is always increased by one (so if they were pentagons, they are no more small); two new pentagons or one new quadrangle are created.
}
\label{fig:oper1}
\end{figure}

We may suppose that $P$ has two small faces $f_1$ and $f_2$ which are not adjacent to each other. We claim that $f_1$ and $f_2$ are at mutual position $(1,1)$ and the edge connecting them is incident to a quadrangle:

Suppose $f_1$ and $f_2$ are two small faces in mutual position $(c_1,c_2)$ such that  $c_1\ge c_2\ge 1$ and $c_1\ge 2$. Then there exists a new patch $P^\prime$ with the same boundary vector and the same curvature, but with a bigger number of vertices: $P^\prime$
can be found by inserting two pentagons and $c_1+c_2-3$ hexagons along a shortest path joining $f_1$ and $f_2$. (The path is in $P$ due to convexity of $P$.) By applying this operation, the size of $f_1$ and $f_2$ is increased by one; the mutual position of the two new pentagons is $(c_1-1,c_2-1)$, see Figure \ref{fig:oper1} for illustration. The patch $P^\prime$ is indeed a patch of a Barnette graph, since no pair of adjacent quadrangles can be created this way.

Similarly, if $c_1\ge 3$ and $c_2=0$, then there is a sequence $h_1,\dots,h_{c_1-1}$ of hexagons forming a dual path joining $f_1$ and $f_2$. We subdivide the edge between $f_1$ and $h_1$ and the edge between $h_{c_1}$ and $f_2$ once; we subdivide each edge between $h_i$ and $h_{i+1}$ ($1 \le i \le c_1-2$) twice; we join the new vertices in such a way that $h_1$ and $h_{c_1-1}$ are split into a pentagon and a hexagon and that all other hexagons in the sequence are split into two new hexagons. Again, the size of $f_1$ and $f_2$ is increased by one and a new pair of pentagons at mutual position $(c_1-2,1)$ is created, see Figure \ref{fig:oper1} for illustration.

Analogously, if $(c_1,c_2)=(2,0)$, then there is a hexagon $h$ adjacent to both $f_1$ and $f_2$. To obtain $P^\prime$, it suffices to subdivide the two edges $h$ shares with $f_1$ and $f_2$, respectively, and join the two new vertices by a new edge. This way $h$ is split into two pentagons and the size of $f_1$ and $f_2$ is increased by one, see Figure \ref{fig:oper1} for illustration.

Finally, let $(c_1,c_2)=(1,1)$. Then $f_1$ and $f_2$ are connected by an edge $e$. If the edge $e$ is not incident to any quadrangle, then new patch $P^\prime$ can be obtained by replacing $e$ by a quadrangle, see Figure \ref{fig:oper1} for illustration. Since the size of $f_1$ and $f_2$ is increased by one, there can not be two adjacent quadrangles in the patch $P^\prime$.

To conclude, for every pair of non-adjacent small faces of $P$, there is a quadrangle adjacent to both of them, so both of them are pentagons, and so $\mu\ge 4$ and $P$ contains a quadrangle adjacent to two pentagons (which are not adjacent to each other). If $\mu=4$, then $P$ has no other small faces, so it is primitive, a contradiction. If $\mu=5$, then $P$ contains an additional pentagon, which, due to the previous observations, has to be adjacent to (the only) quadrangle -- again we obtain a primitive patch, a contradiction.
\ep

\begin{corollary}
For a given curvature and given boundary vector, a convex patch with maximal number of vertices has to be a primitive one. 
\end{corollary}

\begin{lemma}
Let $P$ be a convex patch of curvature $\mu\le 5$ and boundary vector $\sigma$. Then there exists a unique primitive patch $\bar{P}(\mu,\sigma)$ with the same curvature and the same boundary vector.
\label{le:unique}
\end{lemma}

Proof. The existence is given by the previous lemma. The uniqueness can be proven by induction, by adding/removing rows of hexagons from a patch, or, alternatively, by considering embeddings of patches onto infinite hexagonal cones. We omit the details. \ep

\begin{lemma}
Let $P$ be a convex patch of perimeter $p$ and curvature $\mu\le 5$. Then $P$ has at most $\frac{p^2}{6-\mu}$ vertices.
\label{lemma:finite}
\end{lemma}

Proof. It suffices to count the numbers of vertices of primitive convex patches. We omit the details. \ep

It is worth mentioning that the bound from Lemma \ref{lemma:finite} is tight only if $\mu(P)\le 2$ and the patch contains at most one small face.

\begin{corollary}
Let $P$ be a convex patch of curvature $\mu\ge 7$. Then $P$ can be realized only in finitely many Barnette graphs. 
\label{cor:howto}
\end{corollary}

The largest Barnette graph containing a given realizable convex patch of curvature $\mu \ge 7$ can be found by adding the corresponding (unique) primitive patch of curvature $12-\mu$.

\subsubsection{Patch closure and essential patches}

A \emph{$k$-disc} centered at a face $f$ of a plane graph $G$, denoted by $B_k(f)$, is a subgraph of $G$ composed of facial cycles of faces at (dual) distance at most $k$ from the face $f$. Note that if $k$ is large enough, then $B_k(f)=G$ for any $f$.

A patch $P^\prime$ is called a \emph{closure} of a patch $P$, if 
\begin{enumerate}
\item $P$ is contained in $P^\prime$,
\item every small face of $P'$ corresponds to a small face of $P$,
\item $P^\prime$ contains the 2-discs centered at the small faces of $P$, and 
\item $P'$ is convex. 
\end{enumerate}

A patch $P$ is called \emph{closed} if it is a closure of itself.

Clearly, if $P^\prime$ is a closure of $P$, then $P^\prime$ can be obtained from $P$ by adding a finite number of hexagons.

Let $P$ be a patch with boundary vector $\sigma(P)=\sigma_1\sigma_2\dots \sigma_k$. 
The \emph{small face distance} of a value $\sigma_i$ is equal to the minimum of the distances $d(f^*,g^*)$, where $f$ is the new face of the patch $P^{i\gets j}$ (for some $j$ sufficiently big),  $g$ run the set of small faces of $P$, and the distances are taken in the inner dual (dual without the vertex representing the outer face) of $P^{i\gets j}$.

Let $P$ be a patch which is not convex. Then we set all the values of its boundary vector as \emph{admissible}.

Let $P$ be a convex patch with boundary vector $\sigma(P)=\sigma_1\sigma_2\dots \sigma_k$. 
A value $\sigma_i$ is called \emph{admissible}, if the small face distance of $\sigma_i$ is at most 2. 

Observe that boundary vectors of closed patches have no admissible values.

Let $P$ be a patch with boundary vector $\sigma(P)=\sigma_1\sigma_2\dots \sigma_k$ which is not closed. 
A \emph{critical element} of the boundary vector of $P$ is an admissible value $\sigma_i$ such that
\begin{itemize}
\item $\sigma_i$ is maximal, and then
\item the sum $\sigma_{i-1}+\sigma_{i+1}$ (incides taken modulo $k$) is maximal, unless $\mu(P)\ge 5$ and $\max_{i=1}^k\sigma_i=3$; in which case we choose $\sigma_i=3$ contained in a subsequence $32^j3$ of minimum length, and then
\item the small face distance of $\sigma_i$ is minimal.
\end{itemize}

\begin{lemma}
Let $G$ be a Barnette graph and let $f$ be a small face of $G$. Then there exists a finite sequence of patches $\{P_k\}_{k=1}^t$ contained in $G$ such that
\begin{itemize}
\item $P_1$ is a cycle of length equal to the size of $f$;
\item $P_{k+1}=P_k^{i\gets j}$, where $j\in\{4,5,6\}$ and $\sigma_i$ is a critical element of $\sigma(P_k)$;
\item for each $k$, in the embedding of $P_k$ into $G$, the face $f$ corresponds to a face of $P_k$;
\item either $P_t$ is the first closed patch of the sequence or $P_t=G$.
\end{itemize}
\label{le:seq}
\end{lemma}

Proof. The existence of the sequence is guaranteed by Lemma \ref{le:real}. Either adding faces one by one yields a closed patch, or all the faces of $G$ are eventually added. In both cases the sequence is finite. \ep

Observe that the sequence $\{P_k\}_{k=1}^t$ of patches contained in a Barnette graph $G$ starting with a fixed small face $f$ of $G$ given by Lemma \ref{le:real} is not unique -- it may depend on the choice of a critical element.

Let $f$ be a small face of a Barnette graph $G$ and let $P$ be a patch. If $P=P_t$ for some sequence described in Lemma \ref{le:real} starting with $f$, then we call $P$ an \emph{essential patch} for $f$ in $G$. 

\subsection{Patches and the general procedure}

Let $P$ be a patch essential from some small face of a Barnette graph $G$.
Then the Hamiltonian cycle $C_T$ of the triangulation $T$ capturing the mutual position of small faces of $G$ enters and leaves $P$ at least once. 

We will modify the general procedure in order to ensure that we can choose a cycle $C_T$ enterling and leaving $P$ exactly once: For essential patches of curvature at least 6 this is automatically true due to convexity of the patch and minimality of the cycle. For each essential patch $P$ of curvature at most 5 we can temporarily replace $P$ by the corresponding primitive patch $\bar{P}$; in the resulting graph we find the cycle $C^*$ visiting each small face exactly once. Since in $\bar{P}$ the small faces are adjacent to each other, they are consecutive along $C^*$ by minimality of $C^*$.
When replacing back the primitive patches by the actual patches, we keep the order in which the (primitive) patches were covered by $C^*$ and we keep the position of the segments joining different patches. We disregard the way how $C^*$ visits the small faces inside each essential patch, since we will inspect that in details later.

From this point on we may assume that for each essential patch $P$ there are exactly two segments leaving $P$, say $P_i^*$ and $P_j^*$. For any position of the segments $P_{i+1}^*, \dots, P_{j-1}^*$ inside $P$, the difference $\varphi_j \circ \varphi_i^{-1}$ is a permutation of three elements which is even if and only if $\mu(P)$ is even (each pentagon of $P$ contributes with a single transposition). 

If $\mu(P)$ is odd, then the difference $\varphi_j \circ \varphi_i^{-1}$ is an odd permutation -- a transposition. Therefore, among the nine choices of colorings of $G_1$ and $G_2$, for one choice both segments leaving $P$ are inactive, for four choices one of them is active and the other one is inactive, and for the remaining four both segments are active -- the patch behaves like a pentagon. We will call these patches type 1.

If $\mu(P)$ is even and the difference $\varphi_j \circ \varphi_i^{-1}$ is the identity, then among the nine choices of colorings of $G_1$ and $G_2$, for three of them both segments are inactive and for the remaining six both segments are active -- the patch behaves like a quadrangle. We will call these patches type 0.

If $\mu(P)$ is even and the difference $\varphi_j \circ \varphi_i^{-1}$ is an even permutation different from the identity, then it has to be a cycle of length three. Therefore, among the nine choices of colorings of $G_1$ and $G_2$, for three of them both segments are active and for the remaining six there is one active and one inactive segment -- the patch behaves like a pair of pentagons of different colors. We will call these patches type 2.

Let $P$ be a patch essential from some small face of a Barnette graph $G$.
Let the position of two segments leaving $P$ and all the segments inside $P$ be fixed.
Let one of the nine colorings of $G_1$ and $G_2$ be chosen.
 Let the procedure described in Section \ref{sec:pro} be applied. We first obtain a $2^*$-factor, which is then transformed into at most two $2$-factors (depending on the order of decisions at 2-cycles of the $2^*$-factor). 

Let $H^0_P$ be the subgraph of the residual graph $H$ induced by the vertices corresponding to the faces of $P$ and faces adjacent to faces of $P$ in $G$. There can be vertices of degree 1 or 2 in $H^0_P$. We add $3-d$ new vertices adjacent to each vertex of degree $d$ in $H^0_P$ inside the outer face; we then connect all these new vertices by a new cycle. This way we obtain a plane cubic graph $H_P$, we call it \emph{partial residual graph}. 

Let $f^*$ be a vertex of $H_P$ corresponding to a face $f$ of $P$. The face $f$ is a white resonant hexagon. If we recolor $f$ black, then three different components of the underlying $2$-factor are merged into a single cycle; the vertex $f^*$ is deleted from $H_P$ and the three resulting $2$-vertices are suppressed. We call this operation \emph{elimination} of $f^*$.

We say that a plane cubic graph is \emph{strongly essentially $4$-edge-connected}, if it is cyclically 3-edge-connected, and every cyclic 3-edge-cut separates a triangle adjacent to the outer face from the rest of the graph.

We say that a plane cubic graph is \emph{essentially $4$-edge-connected} if it can be transformed into a strongly essentially $4$-edge-connected plane graph by a vertex elimination.

\begin{figure}[pht]
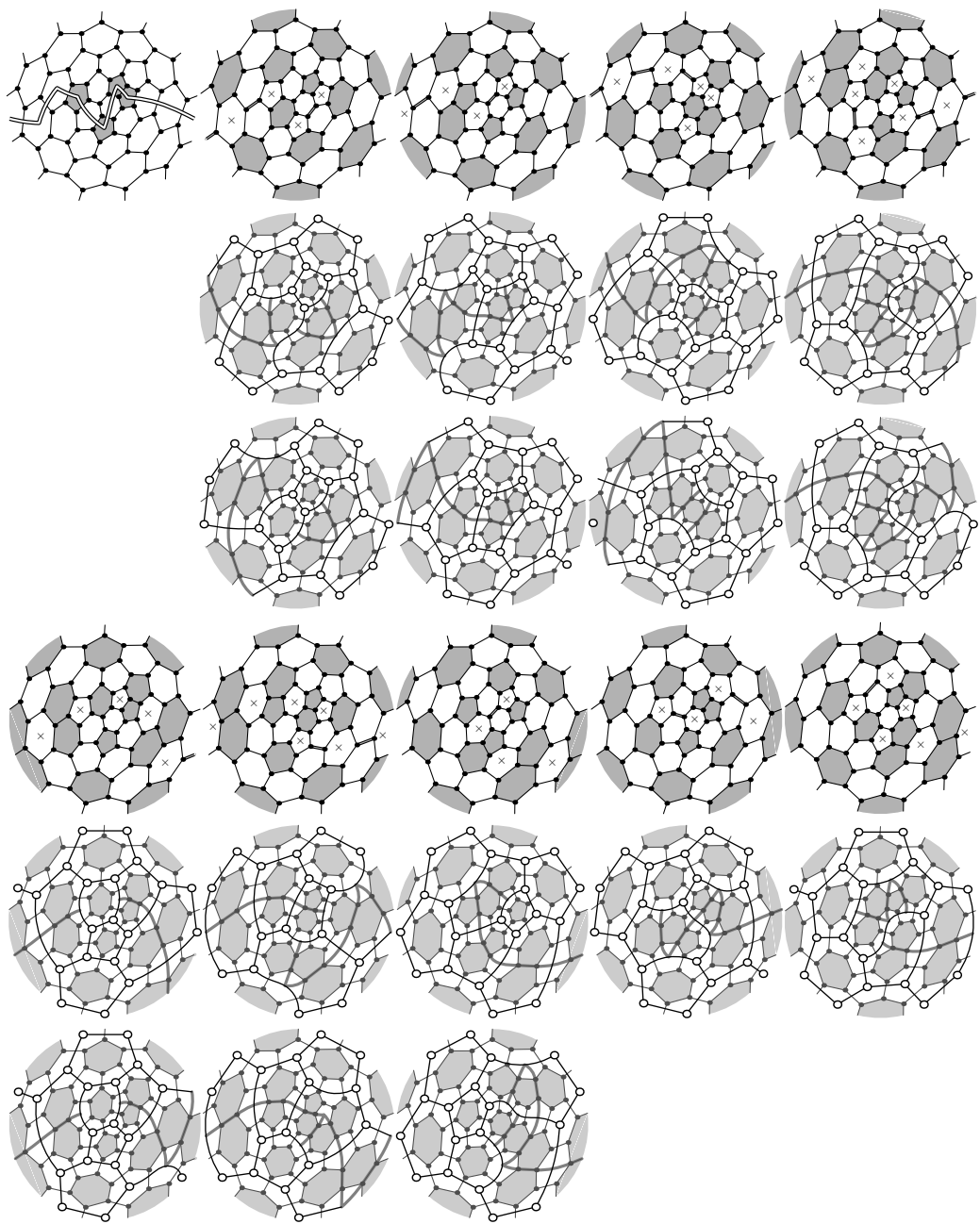

\centerline{
\begin{tabular}{c@{}c@{}c@{}c@{}c}
{\includegraphics[scale=.5]{0.4.24.0}}
&
 {\includegraphics[scale=.5]{0.4.24.2.25_41.32}}
&
 {\includegraphics[scale=.5]{0.4.24.2.36_41.13}}
&
 {\includegraphics[scale=.5]{0.4.24.2.26_41.21}}
&
 {\includegraphics[scale=.5]{0.4.24.2.26_14.12}}
 \\
& {\includegraphics[scale=.5]{0.4.24.2.25_41.320}}
& {\includegraphics[scale=.5]{0.4.24.2.36_41.130}}
& {\includegraphics[scale=.5]{0.4.24.2.26_41.210}}
& {\includegraphics[scale=.5]{0.4.24.2.26_14.120}}
 \\
& {\includegraphics[scale=.5]{0.4.24.2.25_41.321}}
& {\includegraphics[scale=.5]{0.4.24.2.36_41.131}}
& {\includegraphics[scale=.5]{0.4.24.2.26_41.211}}
& {\includegraphics[scale=.5]{0.4.24.2.26_14.121}}
\end{tabular}
}

\centerline{
\begin{tabular}{c@{}c@{}c@{}c@{}c}
 {\includegraphics[scale=.5]{0.4.24.2.25_27.23}}
&
 {\includegraphics[scale=.5]{0.4.24.2.36_30.31}}
& 
 {\includegraphics[scale=.5]{0.4.24.2.27_41.13}}
&
 {\includegraphics[scale=.5]{0.4.24.2.14_41.32}}
&
 {\includegraphics[scale=.5]{0.4.24.2.32_41.21}}
 \\
 {\includegraphics[scale=.5]{0.4.24.2.25_27.231}}
& {\includegraphics[scale=.5]{0.4.24.2.36_30.310}}
& {\includegraphics[scale=.5]{0.4.24.2.27_41.130}}
& {\includegraphics[scale=.5]{0.4.24.2.14_41.320}}
& {\includegraphics[scale=.5]{0.4.24.2.32_41.210}}
 \\
 {\includegraphics[scale=.5]{0.4.24.2.25_27.230}}
& {\includegraphics[scale=.5]{0.4.24.2.36_30.311}}
& {\includegraphics[scale=.5]{0.4.24.2.27_41.139}}
&& 
\end{tabular}
}
\caption{For a patch $P$ with four pentagons and a fixed position of two segments leaving $P$, for each of the nine colorings of $G_1$ and $G_2$ the $2^*$-factor and (at most) two simple $2$-factors obtained by the general procedure are depicted. The third drawing in the third column of the second row proves that for this position of the segments leaving $P$ the patch $P$ is parity-switching.
}
\label{fig:reg}
\end{figure}

\begin{figure}[pht]
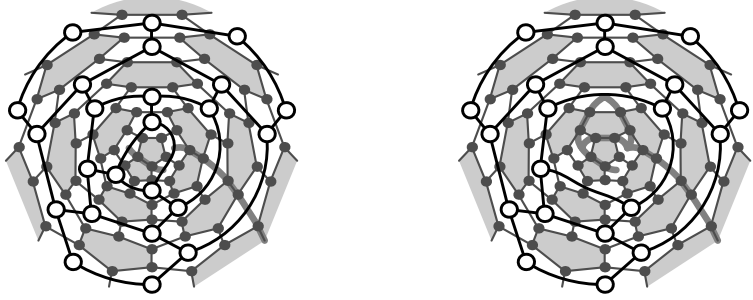

\centerline{
\includegraphics{0.5.prvy.138}
\hfil
\includegraphics{0.5.prvy.38}
}
\caption{For a few patches with many small faces adjacent to each other, the first outcome of the general procedure is a 2-factor such that the corresponding partial residual graph is not strongly essentially 4-edge-connected (left). However, to obtain a strongly essentially 4-edge-connected graph, it suffices to eliminate a vertex incident to a short cycle (right).
}
\label{fig:0.5.prvy}
\end{figure}

We say that a patch $P$ is \emph{regular}, if for every possible position of a pair of segments leaving $P$ and for every choice of the colors of $G_1$ and $G_2$, there is a permutation of small faces of $P$ such that for each of the (at most) two $2$-factors obtained by the general procedure the corresponding partial residual graph is essentially 4-edge-connected. See Figures \ref{fig:reg} and \ref{fig:0.5.prvy} for illustration.

We say that a patch $P$ is \emph{weakly regular}, if for every possible position of a pair of segments leaving $P$ there exists a choice of the colors of $G_1$ and $G_2$ such that there is a permutation of small faces of $P$ such that for at least one $2$-factor obtained by the general procedure the corresponding partial residual graph is essentially 4-edge-connected.

We say that a patch $P$ is \emph{parity-switching} if for every possible position of a pair of segments leaving $P$ there exists a choice of the colors of $G_1$ and $G_2$ such that there exists a permutation of small faces of $P$ such that one of the operations $O_1$ and $O_2$ can be applied inside $P$; for both $2^*$-factors (before and after the operation), for at least one $2$-factor the corresponding partial residual graph is essentially 4-edge-connected.

\subsection{Generation of patches}

\begin{theorem}
There exists a finite set $\mathcal{P}$ of patches such that for every Barnette graph $G$ on at least 318 vertices and every small face $f$ of $G$, there exists a patch $P\in \mathcal{P}$ essential for $f$ in $G$.
\label{th:gen}
\end{theorem}

Proof. We prove the claim by construction. We used Algorithm \ref{algo:add} to generate all the patches in $\mathcal{P}$, by two calls of the procedure \textsc{Generate()}, passing as a parameter first a $4$-cycle and then a $5$-cycle, with the database of patches containing initially the closures of the two initial patches. The procedure uses Algorithm \ref{algo:closure} as a subroutine to calculate a closure of a given patch. 

If the insertion at lines $6$, $11$, or $16$ of Algorithm \ref{algo:add} fails, it means that there is no Barnette graph containing the current patch $P$ such that the element $\sigma_i$ corresponds to a $j$-face for $j=4$, $5$, or $6$, respectively. If this is the case, the following lines are ignored until the next insertion.
Similarly for the insertion at line 6 of Algorithm \ref{algo:closure}.
\ep

\begin{algorithm}[pht]
\caption{Generation of all closed patches containing a given patch}\label{algo:add}
\begin{algorithmic}[1]
\Procedure{Generate}{patch $P$}
\If{$\mu(P)\ge 7$ and the largest graph containing $P$ has at most 316 vertices}
  \Return{}  
\Else 
  \State let $\sigma_i$ be a critical element of the boundary of $P$
  \If {the path along $\sigma_i$ is not adjacent to a 4-face}
    \State $P^\prime \gets P^{i\gets 4}$
    \State $P^{\prime\prime} \gets $\Call{Closure}{$P^\prime$}
    \If{$P^{\prime\prime}$ is not in the database of patches}
  	  \State Add $P^{\prime\prime}$ to the database of patches
  	  \State \Call{Generate}{$P^\prime$}
    \EndIf
  \EndIf  
  \State $P^\prime \gets P^{i\gets 5}$ 
  \State $P^{\prime\prime} \gets $\Call{Closure}{$P^\prime$}
  \If{$P^{\prime\prime}$ is not in the database of patches}
  	\State Add $P^{\prime\prime}$ to the database of patches
  	\State \Call{Generate}{$P^\prime$}
  \EndIf
  \State $P^\prime \gets P^{i\gets 6}$ 
  \State $P^{\prime\prime} \gets $\Call{Closure}{$P^\prime$}
  \If{$P^{\prime\prime}\ne
  P^\prime
  $
  }
	  \State \Call{Generate}{$P^\prime$}
  \EndIf	  
\EndIf
\EndProcedure
\end{algorithmic}
\end{algorithm}

\begin{algorithm}[pht]
\caption{Computation of a closure of a given patch}\label{algo:closure}
\begin{algorithmic}[1]
\Procedure{Closure}{patch $P$}
\If{$P$ is closed}
	  \State \Return{$P$}  
\Else{}
  \State let $\sigma_i$ be a critical element of the boundary of $P$
  \State $P^\prime \gets P^{i\gets 6}$ \label{lineno:add}
  \State \Return{\Call{Closure}{$P^\prime$}}
\EndIf
\EndProcedure
\end{algorithmic}
\end{algorithm}

The counts of patches generated in the proof of Theorem \ref{th:gen} are depicted in Table \ref{table}.

\begin{table}[pht]
\centerline{
\begin{tabular}{|c||c|c|c|c|c|c|c|}\hline
$f_4$ $\backslash$ $f_5$&1&2&3&4&5&6&7 \\\hline\hline
0&1&3&12&92&1202&8821&679 \\\hline
1&3&24&354&3279&&& \\\hline
2&37&383&&&&& \\\hline
\end{tabular}}
\caption{Numbers of essential patches in $\mathcal{P}$, given number of pentagons and quadrangles. Amongst the patches of curvature greater than 6, only patches contained in at least one graph on at least 318 vertices are counted.}
\label{table}
\end{table}

\subsection{Analyse of patches}

The following statements were checked by computer:

\begin{theorem}
There is no patch $P\in \mathcal{P}$ with $\mu(P)\ge 8$. 
\label{th:8+}
\end{theorem}

This means that for every patch $P$ of curvature at least 8 considered by the generating algorithm, the largest graph containing $P$ has less than 318 vertices. See line 2 of Algorithm \ref{algo:add} and the remark after Corollary \ref{cor:howto}.

\begin{theorem}
Every patch $P\in \mathcal{P}$ with $\mu(P)\le 5$ is regular. Every patch $P\in \mathcal{P}$ with $\mu(P)\in \{6, 7\}$ is weakly regular, unless $P$ contains a cap of a nanotube of type $(p_1,p_2)$ with $(p_1,p_2)\in \{(4,0), (5,0),(4,1),(5,1),(3,2),(4,2),(3,3),(4,3)\}$.
\label{th:67}
\end{theorem}

Theorem \ref{th:67} guarantees the existence of a simple $2$-factor such that the residual graph is cyclically 4-edge-connected. The only missing part is that we cannot be sure that this $2$-factor is odd.

We do not need to check for regularity of patches of curvature 6 and 7, weak regularity suffices instead: If a Barnette graph contains an essential patch of curvature $\mu=7$, then it only contains one. Therefore, we can chose the coloring of $G_1$ and $G_2$ such that no segment leaving $P$ is active. 

If a Barnette graph $G$ contains an essential patch $P$ of curvature $\mu=6$, then $P$ contains a cap of a nanotube, and we can choose the coloring of $G_1$ and $G_2$ such that the tubical part of $G$ is traversed by at most one active segment (one if $P$ is type 2, none if $P$ is type 0).

If $G$ is a nanotube of type $(c_1,c_2)$, $c_1\ge c_2\ge 0$, and the caps are (contained in) patches of type 0, then $3 \mid (c_1-c_2)$, so we can write $(c_1,c_2)=(3a+b,b)$ for some integers $a,b\ge 0$. If we choose any of the three colorings of $G_1$ and $G_2$ such that no active segment traverses the tubical part of $G$, then the residual graph $H$ is a nanotube of type $(a+b,a)$.

If $G$ is a nanotube of type $(c_1,c_2)$, $c_1\ge c_2\ge 0$, and the caps are (contained in) patches of type 2, then $3 \nmid (c_1-c_2)$, so we can write $(c_1,c_2)=(3a+b,b+1)$ or $(c_1,c_2)=(3a+b+1,b)$ for some integers $a,b\ge 0$, or $(c_1,c_2)=(3a+2,0)$ for some integer $a\ge 1$. If we choose a coloring of $G_1$ and $G_2$ such that one active segment traverses the tubical part of $G$, then the residual graph $H$ is a nanotube of type $(a+b,a)$ (in the first two cases) or $(a+1,1)$ (in the third case).

This is the reason for excluding the aforementioned 8 types of nanotubes.

\begin{figure}
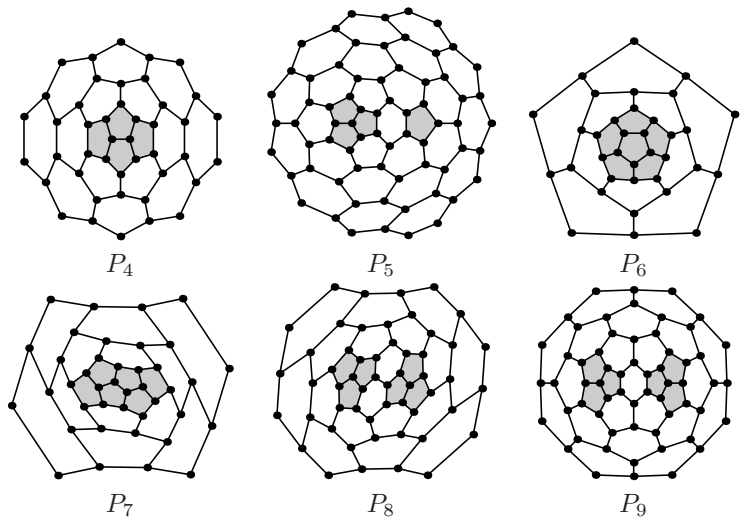

\centerline{
\begin{tabular}{ccc}
\includegraphics[scale=.75]{clusters.4}
&
\includegraphics[scale=.75]{clusters.5}
&\includegraphics[scale=.75]{clusters.6}
\\
$P_4$ & $P_5$ & $P_6$ \\
\includegraphics[scale=.75]{clusters.7}
&
\includegraphics[scale=.75]{clusters.8}
&
\includegraphics[scale=.75]{clusters.9}
\\
$P_7$ & $P_8$ & $P_9$ 
\end{tabular}
}
\caption{Patches with 4 and 6 pentagons for which it is not possible to increase or decrease the number of black pentagons by 2.}
\label{fig:clusters}
\end{figure}

\begin{theorem}
Let $P\in \mathcal{P}$ with $\mu(P)\le 7$. Then $P$ is parity-switching, unless 
$P$ is one of the following exceptional patches:
\begin{itemize}
\item the patch $P_1$ of curvature $1$ having one pentagon,
\item the patch $P_2$ of curvature $2$ containing two adjacent pentagons,
\item the patch $P_3$ with three pentagons sharing a common vertex,
\item two patches $P_4$ and $P_5$ with four pentagons (the type 0 patches obtained from $P_3$ by adding a pentagon at distance at most two),
\item four patches $P_6$, $P_7$, $P_8$, $P_9$ with six pentagons, depicted in Figure \ref{fig:clusters}.
\end{itemize}
\label{th:psw}
\end{theorem}

There is a combinatorial reason for the patches $P_3$-$P_9$ not to be parity-switching: if three pentagons share a vertex, either one or two of them have to be black, so we do not have the freedom to change their colors independentely.

As a consequence of Theorem \ref{th:psw}, if a Barnette graph contains at least one parity-switching essential patch, we choose the coloring of $G_1$ and $G_2$ that allows to change the parity of the $2$-factor, and, by regularity, we are done.

It remains to consider Barnette graphs (in fact, fullerene graphs) only containing patches $P_1$-$P_9$ and verify that we can use parity-switching operations using pentagons from different patches.

If a fullerene graph contains $P_6$, it is a nanotube of type $(5,0)$, and it is known to be Hamiltonian \cite{mar2}. If a fullerene graph contains $P_7$ ($P_8$, $P_9$, respectively), then it is a nanotube of type $(4,2)$ (of type $(6,2)$, $(8,0)$) -- the patch itself already contains a corresponding ring. Out of all the possible patches (caps) to close the other end of the tube, $P_7$ ($P_8$, $P_9$) is the only one that is not parity-switching, as it was verified by a computer. However, if both caps of a nanotube are $P_7$ ($P_8$, $P_9$), then it has an even number of hexagons and exactly 6 black and 6 white pentagons, so by Lemma \ref{l:paths} the number of cycles in the $2$-factor is odd.

It remains to consider fullerene graphs only having patches $P_1$, $P_2$, $P_3$, $P_4$, and $P_5$. 

It was verified by computer that for each of the five patches, for each active segment leaving the patch, the $\Gamma$-path can be transformed into a pair of $\Gamma$-paths (interconnected inside the patch or not) -- it is nothing else than applying a half of one of the operations $O_1$ and $O_2$ (or its inverse) inside the patch and the other half inside another.

In each of the patches this modification corresponds to increasing or decreasing the number of black pentagons by one. In most of the cases both are possible. More precisely, for each segment leaving $P_1$, $P_2$, $P_4$ or $P_5$, out of the nine possible colorings, for three colorings the segment is inactive, for at least two colorings it is possible to increase the number of black pentagons by one, and for at least four colorings it is possible to decrease the number of black pentagons by one. It means that if two of these patches are consecutive along $C^*$, then there exists a coloring such that we can decrease the number of black pentagons in each of them by one.

On the other hand, for $P_3$, it is possible to increase the number of black pentagons by one for four colorings and decrease it for two of them. Again, if there are two such patches consecutive along $C^*$, there exists a coloring such that we can increase the number of black pentagons in each of them by one.

It remains to consider fullerene graphs such that along $C^*$, the patches $P_3$ alternate with other types of patches among $\{P_1,P_2,P_4,P_5\}$. Since each $P_3$ contains three pentagons and there are twelve pentagons altogether, it is easy to see that the number of $P_3$ patches is either 2 or 3. 

If there are two $P_3$ patches, the other two patches have six pentagons, and hence one of them is $P_2$ and the other one is either $P_4$ or $P_5$. 
The patch with four pentagons has to be far from each of the $P_3$ patches, otherwise the graph would be too small (see Lemma \ref{lemma:finite}). The patches $P_4$ and $P_5$ are both type 0. That is why we may omit the four-pentagon patch and search only for a cycle passing through the eight pentagons of the other three patches; we consider $P_4$ or $P_5$ as if no segment leaving it was active. As a consequence, we find two $P_3$ patches consecutive along $C^*$.

If there are three $P_3$ patches, the other three patches can only have one pentagon each. 
Moreover, the condition that for each segment joining a $P_3$ to a $P_1$ the two colorings allowing to decrease the number of black pentagons in $P_3$ correspond to the two colorings allowing to increase the number of black pentagons in the other patch implies that out of the nine colorings, there is one with no active segment joining a $P_3$ to a $P_1$, there are four colorings with three active segments and three inactive segments alternating, and there are four colorings with all the six segments active. In all the cases there are three $\Gamma$-paths in $G$. (In the case of no active segments joining different patches, there is still a $\Gamma$-path joining different pentagons inside each $P_3$.)

If we replace a vertex incident to three pentagons inside each $P_3$ by a triangle temporarily, then the graph will contain three pentagons and three triangles (and all the other faces will be hexagons). Moreover, in the coloring of $G_1$ and $G_2$ all the six small faces have the same color.

By the structural theorem of Alexandrov, such a graph can be isometrically embedded onto a surface of a (possibly degenerate) convex polyhedron, say $P$. The polyhedron $P$ has six vertices, and the cycle $C^*$ is a Hamiltonian cycle in some triangulation of $P$.	

\begin{figure}[pht]
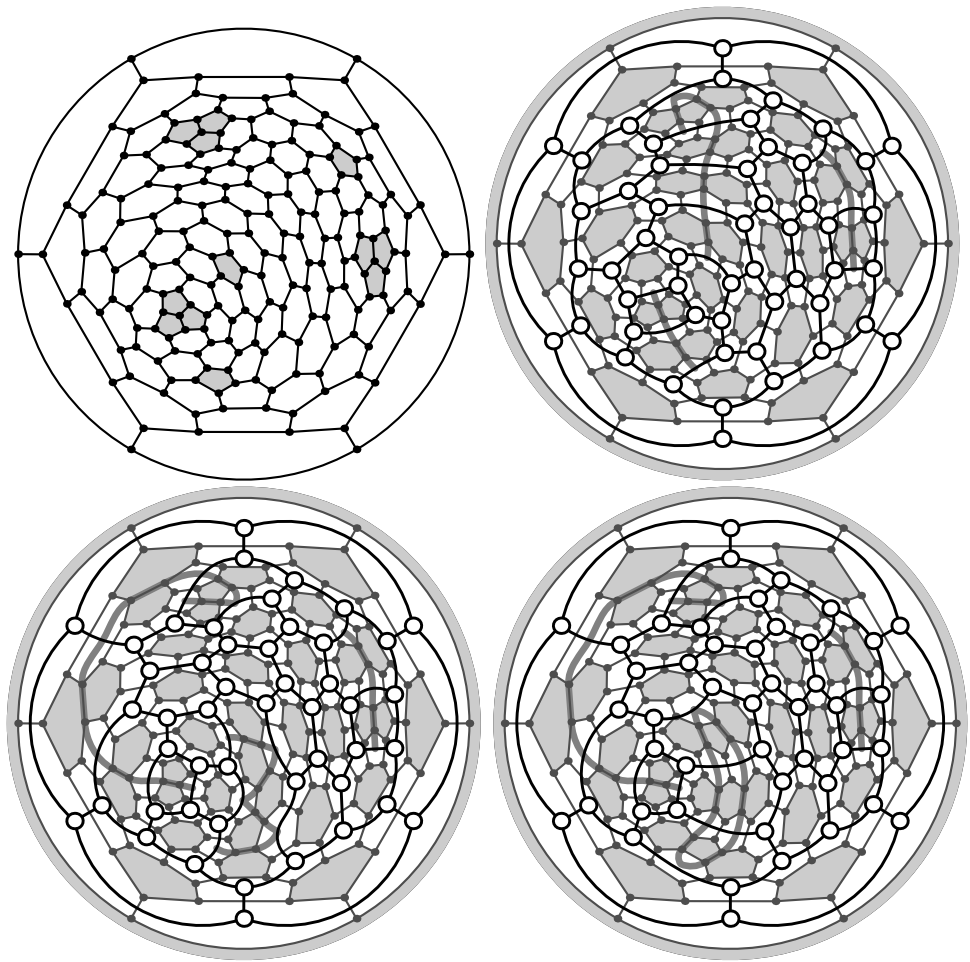

\centerline{
\includegraphics[scale=1]{313131.0}
\includegraphics[scale=1]{313131.1}
}

\centerline{
\includegraphics[scale=1]{313131.2}
\includegraphics[scale=1]{313131.3}
}
\caption{Top to bottom, left to right: An exemple of a fullerene graph on 198 vertices containing three patches $P_3$ and three patches $P_1$. An even 2-factor with three $\Gamma$-paths without a possibility to apply $O_1$ or $O_2$. Another even 2-factor obtained by switching the order of the patches. An odd 2-factor after applying $O_1$.}
\label{fig:313131}
\end{figure}

The cycle $C^*$ cuts the polyhedron $P$ into two hexagons. In the two hexagons the angles at a fixed $P_3$-vertex (center of a triangle) sum up to $180^\circ$, and hence they are both always convex (smaller than $180^\circ$). For the angles at the $P_1$-vertices (centers of isolated pentagons), in at least one hexagon the angle is convex. Therefore, it is always possible to permute a $P_1$ patch with a $P_3$ patch to obtain a new cycle with two consecutive patches of the same type, which gives us a possibility to change the parity of the number of cycles. See Figure \ref{fig:313131} for illustration.

\section{Concluding remarks}

Similar technique could be used to prove Hamiltonicity of related graph classes: planar cubic graphs with only a few faces of size larger than six; projective-planar graphs with faces of size at most six (except, of course, for the Petersen graph), etc.

\end{document}